%% file: Manuscript.tex
\documentclass[3p, number,sort&compress,times,fleqn]{elsarticle}

\usepackage{amsmath}
\usepackage{amssymb}
\usepackage{mathtools}
\usepackage[colorlinks = true,
            linkcolor = blue,
            urlcolor  = blue,
            citecolor = blue,
            anchorcolor = blue]{hyperref}

\usepackage[utf8]{inputenc}
\usepackage{subfig}

\usepackage{csml}

\begin{document}
	
	\begin{frontmatter}
		
		\title{Total Lagrangian Smoothed Particle Hydrodynamics with An Improved Bond-Based Deformation Gradient for Large Strain Solid Dynamics}
		\author[1]{I. M. Wiragunarsa}
		\author[1]{L. R. Zuhal \corref{cor1}}
  		\ead{lavi.zuhal@itb.ac.id}
		\author[1]{T. Dirgantara}
		\author[1]{I. S. Putra}
		\author[1,2]{E. Febrianto}

		\cortext[cor1]{Corresponding author}
		
	   \address[1]{Faculty of Mechanical and Aerospace Engineering, Institut Teknologi Bandung, Ganesa Street 10, Bandung, 40132, Indonesia}
	   \address[2]{Glasgow Computational Engineering Centre, University of Glasgow, University Avenue, Glasgow, G12 8QQ, UK}
		
		\begin{abstract}
		Total Lagrangian Smoothed Particle Hydrodynamics (TLSPH) is one variant of SPH where the variables are described using the fixed reference configuration and a Lagrangian smoothing kernel. TLSPH elevates the computational efficiency of the standard SPH when no topological change is involved, and it alleviates the stability of SPH scheme with respect to tensile loading. However, instabilities associated with spurious mode, or hourglass/zero-energy mode, persists and often affects the simulation of solids undergoing extremely large strain. This work proposes an alternative formulation to compute deformation gradient with improved accuracy and therefore minimising the possibility of encountering the zero-energy mode. Specifically, we leverage the local discrete computation of bond-based (or pairwise) deformation gradient smoothed by the kernel. Additionally, the bond of a particle with itself is considered to preserve the polynomial reproducibility imposed by the kernel correction scheme. We showcase the convergence of the approach using a two-dimensional benchmark example. Furthermore, the accuracy, robustness, and stability of the proposed approach are assessed in various two- and three-dimensional examples, highlighting on the stability improvement that allows for solid dynamic simulations with more extreme elongation. 
		\end{abstract}
		
		\begin{keyword}
			SPH \sep Particle Methods \sep Computational Solid Mechanics
		\end{keyword}
		
	\end{frontmatter}
	

\input{introduction}

	\input{discretisation}

	\input{tlsph}

	\input{pseudostaggered}

	\input{examples}

	\input{conclusions}

        \section*{Acknowledgments}
        The authors gratefully acknowledge the Penelitian Dasar administrated 
        by Direktorat Riset dan Pengabdian Masyarakat - Direktorat Jenderal Penguatan Riset dan Pengembangan - Kementerian Riset dan Teknologi/Badan Riset dan Inovasi Nasional, Republik Indonesia, for the support funding of this research. The authors also thank Chun Hean Lee for the stimulating discussion on the topic. 
        
        \appendix

\include{appendix}

	\bibliographystyle{elsarticle-num-names}
	\bibliography{references.bib}
	
\end{document}

%% file: introduction.tex
\section{Introduction}
\label{sec:intro}
\

Smoothed particle hydrodynamics (SPH) is a meshfree method that uses kernel estimates to approximate a continuous field over the domain~\cite{Lucy1977,GingoldMonaghan1977,BonetSUL2007}. In SPH, the continuum is discretised using a finite set of particles, where the discrete values of the fields are sampled. The continuous field is then obtained by accounting the contribution of the neighbouring particles weighted by the kernel. Owing to its meshfree nature, SPH avoids challenges in creating analysis-suitable meshes commonly found in the mesh-based computational methods~\cite{Hughes2005, Liu2010}. Over the years, SPH has been applied to various scientific and engineering problems, including linear and nonlinear solid mechanics~\cite{Libersky1993}, structural dynamics~\cite{ChenDyn2000}, crack propagations~\cite{Ganesh2022, Ganesh2022a}, and problems involving high distortion and impact~\cite{Wingate1993,Gordon1996,Chakraborty2013,Islam2017,Wira2024}. Despite its notable success, SPH is not without its challenges, manifesting in issues such as lack of consistency, tensile instability, and the presence of zero-energy or hourglass mode~\cite{Vignjevic2009,Vacondio2021}. 

The consistency issue in SPH is often attributed to the discrete approximation of the  convolution integral, which prohibits the scheme to reproduce polynomial of order $p$. This limitation becomes more pronounced when the particle distribution is non-uniform, and in the region near the boundaries. To address this challenge, a common strategy in SPH involves introducing corrections to the kernel, ensuring the reproducibility of polynomials up to a certain order~\cite{BonetKulasegaram, frontiere2017}. This approach is also employed in other meshless methods~\cite{liu1995reproducing}. When the governing equation involves only the first derivative, the correction can be applied to the gradient rather than the kernel itself~\cite{BonetVariational1999, Liu2006}. In certain cases, this modification proves to be more computationally efficient and is adequate for achieving linear consistency. However, one drawback of the correction approach is the loss of symmetry in the kernel, and similarly, its gradient is no longer anti-symmetric, impacting the local momentum and energy preservation. To reconcile both consistency and local conservation, one potential solution is to symmetrise the kernel and its gradient, as demonstrated for instance in~\cite{Islam2018}. 

There are two known sources of instability in the classical, or updated Lagrangian, SPH (ULSPH) formulation. The first is instability under tensile loading, and the second is associated with the zero-energy or hourglass mode. Comprehensive stability analyses of ULSPH are presented in~\cite{swegle1995,Belytschko2000}. Various approaches have been proposed to address the issue of tensile instability, including the introduction of artificial viscosity~\cite{MonaghanFreeSurface1994, Owen2004} and artificial stress~\cite{MonaghanTensile2000}, conservative smoothing~\cite{Randles1996,Hicks1997,Hicks2004,Guenther1994}, a stable adaptive kernel~\cite{Lahiri2020}, and the addition of a separate set of points to evaluate derivatives, known as stress points~\cite{Dyka1997}. Despite being the most prevalent in SPH, the artificial viscosity may lead to an excessive dissipation~\cite{Wang2019}. When SPH is formulated with a Lagrangian kernel, the issue of tensile instability can be eliminated without introducing additional terms to the governing equation. This is often referred to as total Lagrangian SPH (TLSPH)~\cite{Bonet2001,Rabczuk2004,Vignjevic2006}. Furthermore, in TLSPH, the field variables are described in the reference configuration, resulting in a computationally efficient scheme as the neighbour search is performed only once. In recent years, TLSPH has found successful applications in various aspects of solid mechanics, including those involving geometric  nonlinearity~\cite{JunLinThickShell2014, JunLinThinShell2014,JunLinGNFGD2018}, solid dynamics \cite{Lee2017,Lee2019}, and simulations of crack growth~\cite{Wira2021}.

While TLSPH effectively resolves the issue of tensile instability, the persistent presence of oscillation due to the zero-energy, or hourglass, modes remains a challenge. The zero-energy mode is observed in various other numerical methods, including finite difference and finite elements~\cite{Flanagan1981,Jacquotte1984}, and is attributed to the inaccuracies in estimating  the gradient of a discretely represented field. Numerous methods have been proposed to address the zero-energy mode, such as $\delta$-SPH diffusive terms~\cite{Antuono2010,Antuono2012,Hammani2020}, the Rieman solver~\cite{Meng2022}, and moving least square stress regularisation~\cite{Nguyen2017,Lian2023}. The most notable treatment for hourglass instability is the stress-point approach~\cite{Vignjevic2000, Randles2000SP}, which determines the field gradient on a separate set of nodes, sharing roots with the staggered grid approach in finite difference~\cite{Armfield1991,Hustedt2004,Liu2009}. 
In the context of TLSPH, a notable approach involves adding a correctional force~\cite{Ganzenmuller2015}, analogous to those proposed in finite elements. A recent work by Wu~\cite{Wu2023} introduces a non-hourglass TLSPH formulation that decomposes the shear deformation and stress according to the volumetric and deviatoric contributions. Despite showing a significant stability improvement with respect to the hourglass mode, both approaches rely on a set of user-prescribed parameters that are problem dependent, which requires thorough parametric studies.

This work proposes an alternative approach to improve the stability of the TLSPH scheme, specifically targeting the spurious zero-energy mode, without depending on user-defined parameters. Our key ingredient is the alternative expression to compute the deformation gradient, which utilises local pair-wise discrete gradient estimates weighted by the Lagrangian kernel. It is noteworthy that the conventional TLSPH approach in computing the deformation gradient involves smoothing the discrete deformation weighted by the gradient of the Lagrangian kernel. Our proposed methodology traces back to the differentiation identity of the convolution of continuous functions, i.e., $D\, (W \ast x) = D \, W \ast x = W \ast D \, x$. A similar concept has also been proposed in peridynamics~\cite{Silling2017, Breitzman2018}, which is referred to as bond-based deformation gradient approach. In contrast to this approach, we include in our proposal an estimate for the deformation gradient at the particle, leveraging the standard gradient estimate in TLSPH. We demonstrate through a simple one-dimensional example that our proposal improves the gradient estimation and preserves the polynomial reproducibility of the scheme, as governed by the corrected kernel. We tailor our proposed deformation gradient estimate in a mixed formulation TLSPH following the framework of~\cite{BonetPart12015, AguirreJST2014, Lee2016}, enhancing it with the Jameson-Schmidt-Turkel (JST) scheme to prevent excessive dissipation in dynamic simulations. Additionally, for time integration, we employ an explicit three-stage time integrator~\cite{Lee2022}, enabling the use of higher CFL numbers and resulting in faster computational times with second-order accuracy~\cite{Jameson2017}. Furthermore, we validate the proposed approach through various benchmark cases in 2D and 3D, showcasing improved accuracy and stability.

The structure of this paper is as follows. First, we revisit the spatial discretisation in SPH including the kernel correction necessary for achieving the desired polynomial reproducibility.  Subsequently, we discuss the total Lagrangian SPH formulation by deriving from the conservation laws and employing the spatial discretisation based on the corrected Lagrangian kernel. In the next section we describe the proposed alternative estimation of the deformation gradient. We provide an illustration of the proposed approach in one-dimensional setting and subsequently detail its integration into our TLSPH framework. Next we present several numerical examples in two and three dimensions, within the context of solid dynamic problems. In the concluding section, we summarise our study on the proposed approach and outline potential future works.

%% file: discretisation.tex
\section{Spatial discretisation in SPH}
\label{sec:discretisation}

In this section, we discuss the spatial discretisation in SPH by considering a scalar field $f(\vec x)$, where $\vec x \in \Omega$ represents spatial coordinates in $\mathbb{R}^d$, with $d \in \{1, 2, 3\}$. In this context, we assume that the domain $\Omega$ remains static and undeformed, eliminating the need to distinguish between the current and reference systems. The discretisation of the field involves employing values obtained at a finite point set, referred to as particles. The relationship between the discrete point-wise values and the continuous form $f(\vec x)$ is established through the SPH discretisation, utilising a kernel with finite support. This section further details the treatments necessary to ensure compliance with the linear polynomial reproducibility conditions. 

\subsection{Field approximation from point values}
\label{sec:field-approx}

We start our discussion with the acquisition of a smooth approximation for an arbitrary field $ f(\vec x)$ through convolutional smoothing with a kernel $W(\vec x)$
\begin{equation}
	\label{eq:conv}
	\widehat{f} (\vec x) = \int_{\Omega} f(\vec x') W(\vec x - \vec x', \, h) \D \vec x' \, .
\end{equation}
The smoothing kernel is typically required to be positive, has a unit volume, and be compactly supported~\cite{Liu2010}. When choosing kernels with compact support, the integration domain in~\eqref{eq:conv} can be simplified from the whole domain $\Omega$ to only the kernel support centred at evaluation point $\vec x$. In this work, we focus on radial kernels with a compact support of radius $R = \kappa \, h$, where $h$ represents the characteristic smoothing length of the kernel. The factor $\kappa$ varies based on the type of kernel used, for example, $\kappa = 2$ for cubic spline and $\kappa=3$ for quintic spline kernel~\cite{Liu2010}. It is noteworthy that as $h$ approaches $0$, the kernel $W$ becomes a Dirac's delta function $\delta(\vec x - \vec x')$. The application of $\delta(\vec x - \vec x')$ as a kernel in~\eqref{eq:conv} yields the identity $\widehat{f} (\vec x) =  f(\vec x)$. However, in general this identity does not hold for arbitrary kernels, and the approximation $\widehat{f}(\vec x)$ deviates from the original function $f(\vec x)$ according to the moments of the kernel $W$. 

\begin{figure}[ht!]
	\centering
	\includegraphics[width=0.35\textwidth]{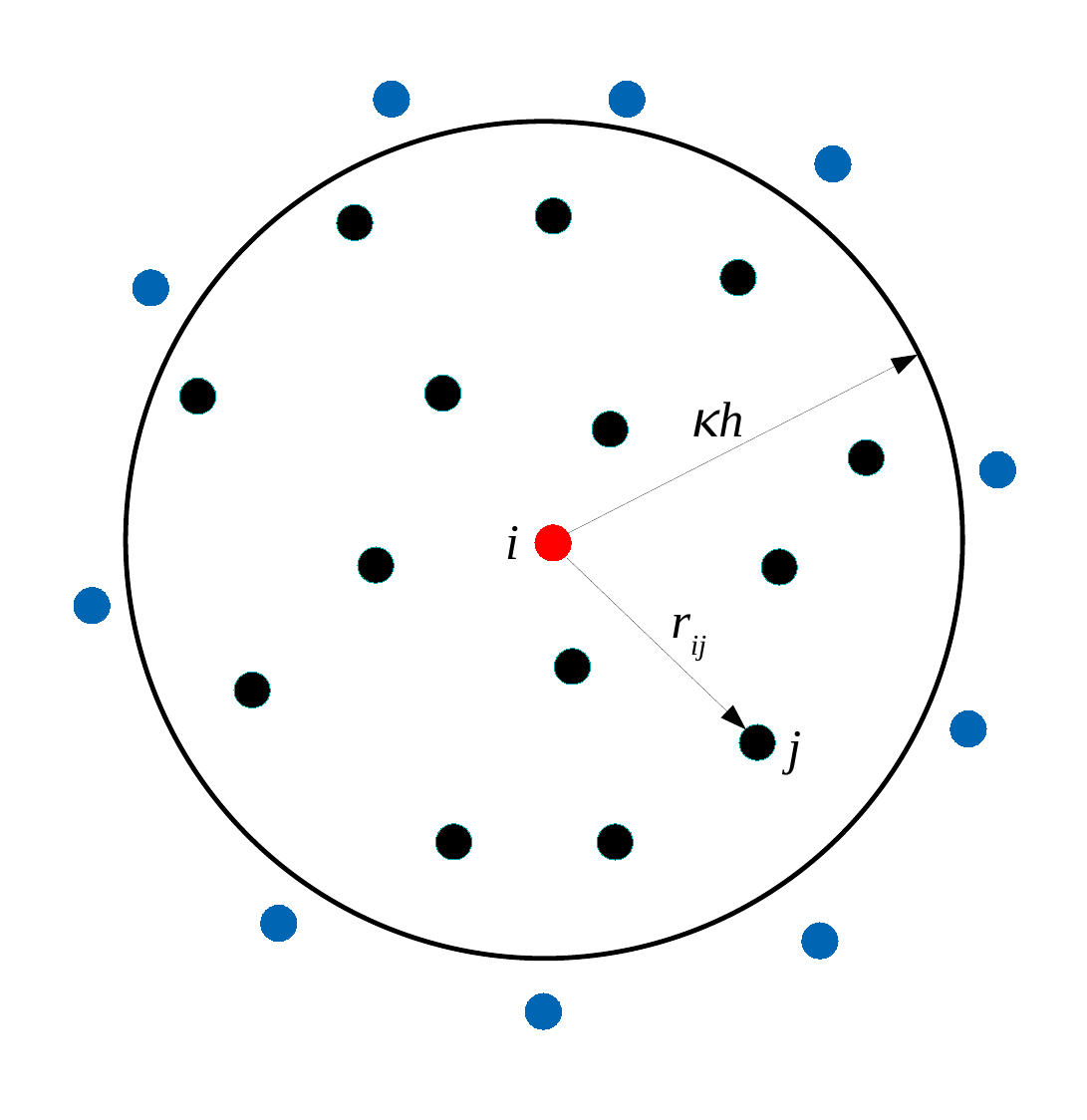}
	\caption{Support of the kernel centred at the $i$-th particle, which is a circle in $\mathbb R^2$ with radius $R = \kappa h$. The neighbouring particles lying inside the support is coloured black.}
	\label{fig:support_domain}
\end{figure}

When the field values are available only at a set of $N$ points $\{\vec x_j\}_{j=1}^N$ scattered across the domain $\Omega$, the convolution integral~\eqref{eq:conv} can be approximated using a summation
\begin{equation}
	\label{eq:sph_f}
	f_h (\vec x) = \sum_{j \in \set N(\vec x)} \frac{m_j}{\rho_j}  f(\vec x_j) W(\vec x-\vec x_j, \, h) \approx  \widehat{ f}(\vec x)  \, ,
\end{equation}
where $m_j$ and $\rho_j$ represent the mass and density of the $j$-th particle, respectively, and their ratio is equivalent to the volume associated with the $j$-th particle, $V_j$. The volume occupied by each particle can be obtained from the Voronoi tessellation of the particle distribution. One possible simplification involves calculating the volume of a spheroid with a radius equals to the average distance between particles. The notation $\set N(\vec x)$ indicates a set of neighbouring particles lying inside the kernel support centred at $\vec x$, see the illustration in Figure~\ref{fig:support_domain}. 

Assuming the kernel is suitably smooth, the gradient of the approximation can be obtained by convolving with the kernel gradient, i.e., 
\begin{equation}
	\label{eq:sph_df}
	\nabla f_h (\vec x) = \sum_{j \in \set N(\vec x)} \frac{m_j} {\rho_j} f(\vec x_j) \nabla W(\vec x - \vec x_j, \, h)	\, .
\end{equation}
Higher-order derivatives can be obtained similarly using higher-order derivatives of the kernel. Due to the radial definition of the kernel, a chain rule is required to obtain the derivative with respect to the coordinate axes. While it is possible to evaluate the approximant~\eqref{eq:sph_f} at any point in the domain $ \vec x \in \Omega $, in SPH, the focus is often on evaluating at a point belonging to the particle set $\vec x_i \in \{ \vec x_j \}$, see Figure~\ref{fig:support_domain}. In this case, the approximation at $\vec x_i$ becomes
\begin{equation}
	\label{eq:sph_disc_i}
	f_h (\vec x_i) = \sum_{j \in \set N(\vec x_i)} V_j  f(\vec x_j) W_j (\vec x_i)  \, . 
\end{equation}
Here, we introduce a simplified notation $W_j(\vec x_i)$ for $W(\vec x_i - \vec x_j)$. The set of neighbouring points within the support of the kernel centred at $\vec x_i$ is indicated as $\set N(\vec x_i)$. These compact notations will be used throughout the remainder of the paper. 

\subsection{Corrected SPH discretisation}
\label{sec:corrected-sph}
The approximation introduced in~\eqref{eq:sph_disc_i} does not inherently ensure the reproduction of constant and linear functions across the domain. To achieve such linear polynomial reproducibility, the following criteria must be met
\begin{subequations}
\begin{align}
	\sum_{j \in \set N(\vec x_i)} V_{j} W_j (\vec x_i) &= 1  \label{eq:zero_complete} \\
	\sum_{j \in \set N(\vec x_i)} V_{j} (\vec x_i - \vec x_j) W_j (\vec x_i) &= 0 \label{eq:first_complete} \, .
\end{align}
\end{subequations}
The first criterion~\eqref{eq:zero_complete} is the discrete manifestation of the kernel's unit volume, which ensures the reproduction of any constant function. Additionally, satisfying both conditions~\eqref{eq:zero_complete} and~\eqref{eq:first_complete} ensures the reproducibility of any linear function. In SPH, it is common practice to modify the smoothing kernel by introducing correction terms and imposing these linear reproducibility criteria. Following~\cite{Liu1995, BonetKulasegaram}, the kernel correction is given by
\begin{equation}
	\label{eq:cw}
	\widetilde{W}_{j}^1 (\vec x_i) = \alpha_i \left( 1 +  \vec{\beta}_i \cdot (\vec x_i - \vec x_j) \right) W_j(\vec x_i) \, .
\end{equation}
Substituting the definition of the corrected kernel~\eqref{eq:cw} into the criterion in~\eqref{eq:first_complete} yields the corrective term $\vec \beta_i$, that is,  
\begin{equation}
	\label{eq:beta}
	\vec{\beta}_i = \left( \sum_{j \in \set N(\vec x_i)} V_{j} (\vec x_i - 
	\vec x_j) \otimes (\vec x_i - \vec x_j) W_j (\vec x_i) \right)^{-1}
	\sum_{j \in \set N(\vec x_i)} V_{j} (\vec x_j - \vec x_i) W_j(\vec x_i) \, .
\end{equation}
Consecutively, the term $\alpha_i$ is obtained by substituting the computed~\eqref{eq:beta} into the definition~\eqref{eq:cw} and the criterion~\eqref{eq:zero_complete}, resulting in
\begin{equation}
	\alpha_i = \frac{1}{\sum_{j \in \set N(\vec x_i)} V_{j} W_j (\vec x_i)
	\left(1 + \vec{\beta}_i \cdot (\vec x_i - \vec x_j) \right) } \, .
	\label{eq:alpha}
\end{equation}
The gradient of the first degree corrected kernel,  $\nabla \widetilde{W}_j^1 (\vec x_i)$, can be obtained explicitly through multiple application of chain rule~\cite{frontiere2017}. The procedure to obtain $\nabla \widetilde{W}_j^1 (\vec x_i)$ is detailed in~\ref{sec:appendix}. The gradient $\nabla \widetilde{W}_j^1 (\vec x_i)$ satisfies the following criteria~\cite{BonetVariational1999}
\begin{subequations}
	\begin{align}
		\sum_{j \in \set N(\vec x_i)} V_{j} \, \nabla W_j (\vec x_i) &= \vec 0 \label{eq:zero_dw} \\
		\sum_{j \in \set N(\vec x_i)} V_{j} \, (\vec x_j - \vec x_i) \otimes \nabla W_j (\vec x_i) &= \vec I  \label{eq:first_dw} \, . 
	\end{align}
\end{subequations}
It is important to note that the corrective terms $\alpha_i$ and $\vec \beta_i$ need to be evaluated at every particle. Moreover, to obtain $\vec \beta_i$, an inversion of small matrix is required at every particle. 

One possible simplification involves modifying the kernel to only satisfy the zero-th order reproducibility~\eqref{eq:zero_complete}, which is equivalent with setting the term $\vec \beta_i = \vec 0$, and subsequently impose an additional criterion on the gradient. This approach leads to a simpler form of the kernel
\begin{equation}
	\widetilde{W}_{j}^0 (\vec x_i)  = \frac{W_j (\vec x_i)}{\sum_{j \in \set N(\vec x_i)} V_{j} W_j (\vec x_i)} \, .
	\label{eq:cw_0}
\end{equation}
The gradient of the zero-th order corrected kernel, $\nabla \widetilde{W}_j^0 (\vec x_i)$, does not satisfy the criterion~\eqref{eq:first_dw}. To address this, a widely adopted approach is to introduce a gradient correction matrix $\vec L_i$ to modify the gradient~\cite{BonetVariational1999}, 
\begin{equation}
	\label{eq:grad_cw}
	\widetilde{\nabla} \widetilde{W}_j^0 (\vec x_i) = \vec L_i \,  \nabla \widetilde{W}_j^0 (\vec x_i) \, ,
\end{equation}
such that criterion~\eqref{eq:first_dw} is satisfied. The correction matrix $\vec L_i$ can be obtained by substituting the above definition~\eqref{eq:grad_cw} into the criterion~\eqref{eq:first_dw}, i.e., 
\begin{equation}
	\label{eq:matL}
	\vec L_i = \left( \sum_{j \in \set N(\vec x_i)} V_{j}(\vec x_j - \vec x_i) \otimes \nabla \widetilde{W}_j^0 (\vec x_i) \right)^{-1} \, .
\end{equation}
The term $\nabla \widetilde{W}_j^0 (\vec x_i)$ in~\eqref{eq:grad_cw} and~\eqref{eq:matL} can be obtained by differentiating~\eqref{eq:cw_0}, resulting in
\begin{equation}
	\nabla \widetilde{W}_{j}^0 (\vec x_i)  = \frac{\nabla W_j (\vec x_i) - \gamma(\vec x_i)}{\sum_{j \in \set N(\vec x_i)} V_{j} W_j (\vec x_i)} \, ,
	\label{eq:grad_cw_0}
\end{equation}
where 
\begin{equation}
	\gamma (\vec x_i)  = \frac{\sum_{j \in \set N(\vec x_i)} V_{j} \nabla W_j (\vec x_i)}{\sum_{j \in \set N(\vec x_i)} V_{j} W_j (\vec x_i)} \, .
	\label{eq:gamma}
\end{equation}
It is noteworthy that the resulted gradient~\eqref{eq:grad_cw_0} satisfies both criteria~\eqref{eq:zero_dw} and~\eqref{eq:first_dw}. This approach is often preferred especially when the discretisation of the governing equations only requires the gradient of the kernel, and not the kernel itself. 

\begin{figure}[ht!]
	\centering
	\includegraphics[width=1\textwidth]{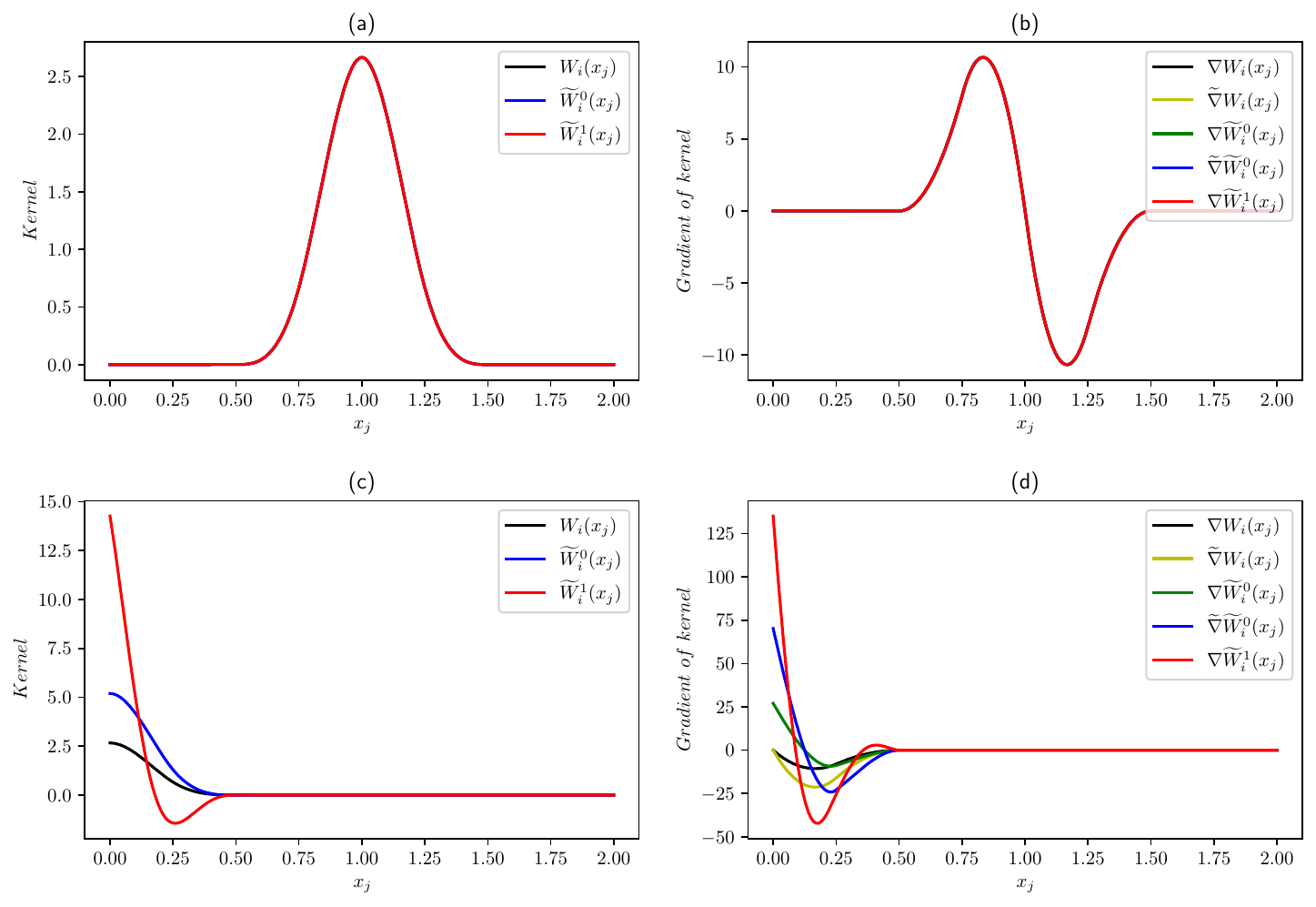}
	\caption{The kernel and kernel gradient of particle $i$ due to interaction with particle $j$: (a) kernel value, and (b) gradient of kernel at $x_i=1$. (c) kernel value, and (d) gradient of kernel at $x_i=0$.}
	\label{fig:kernel}
\end{figure}

As an illustrative example, we visualise the kernel and its gradient in a one-dimensional domain $\Omega = (0, \, 2) \in \mathbb R^1$ with a uniformly distributed set of $N = 200$ particles, considering the different types of correction discussed in this section. For the smoothing kernel $W(x)$, we opt for a cubic B-spline with the support size $\kappa h = 1$. In this context, we compare the original (uncorrected) kernel~$W_j(x_i)$, zero-th order corrected kernel~$\widetilde{W}_{j}^0 (\vec x_i)$, and first order corrected kernel~$\widetilde{W}_{j}^1 (\vec x_i)$, along with their respective first derivatives. Specifically, we investigate the gradient of the original kernel~$\nabla W_j (\vec x_i)$, the corrected gradient of the original kernel~$\widetilde{\nabla} W_j (\vec x_i)$, the gradient of the zero-th order corrected kernel~$\nabla \widetilde{W}_j^0 (\vec x_i)$, the corrected gradient of the zero-th order corrected kernel~$\widetilde{\nabla} \widetilde{W}_j^0 (\vec x_i)$, and finally the gradient of the first order corrected kernel~$\nabla \widetilde{W}_j^1 (\vec x_i)$. 

At $x_i = 1$, Figure ~\ref{fig:kernel}(a) displays the kernels and their gradients are shown in~\ref{fig:kernel}(b). It is observed that in the region far from the boundary, the three kernels overlap, as shown in~Figure ~\ref{fig:kernel}(a). Figure~\ref{fig:kernel}(b) shows that all the gradient variants share the same values and exhibit antisymmetry. However, at the boundary, distinct shapes are evident in the kernels as shown in Figure~\ref{fig:kernel}(c) for $x_i = 0$. Notably, the first order corrected kernel~$\widetilde{W}_{j}^1 (\vec x_i)$ is no longer strictly positive and exhibits a higher peak than the other two variants. Additionally, all three kernels still satisfy the unit volume criterion~\ref{eq:zero_complete} and display asymmetry when $x_i$ is near the boundary. In~\ref{fig:kernel}(d), the gradient~$\nabla \widetilde{W}_j^1 (\vec x_i)$ exhibits the highest peak at $x_i = 0$ and has two inflection points, while all of the gradient variants are no longer antisymmetric. Finally, it is important to note that higher order consistent kernels (and their respective derivatives) generally entail more involved computations. Therefore, in practice, the choice of the correction technique is often justified based on computational cost and the required accuracy.

%% file: tlsph.tex
\section{Total Lagrangian SPH for solid body}
\label{sec:tlsph}

In this section we describe the motion and deformation of an arbitrary solid body and establish the relationship between the current and reference configurations. We then proceed to derive the conservation laws in the reference configuration, closely following the total Lagrangian SPH formulation introduced in~\cite{Lee2016}. Additionally, we discuss in this section the Jameson-Schmidt-Turkel (JST) technique employed to control dissipation in our dynamic simulations, and describe the explicit three-stage time integrator which relaxes the restriction on the time discretisation. 

\subsection{Motion and deformation}
\label{sec:motion-deform}

\begin{figure}[]
	\centering
	\includegraphics[width=0.5\textwidth]{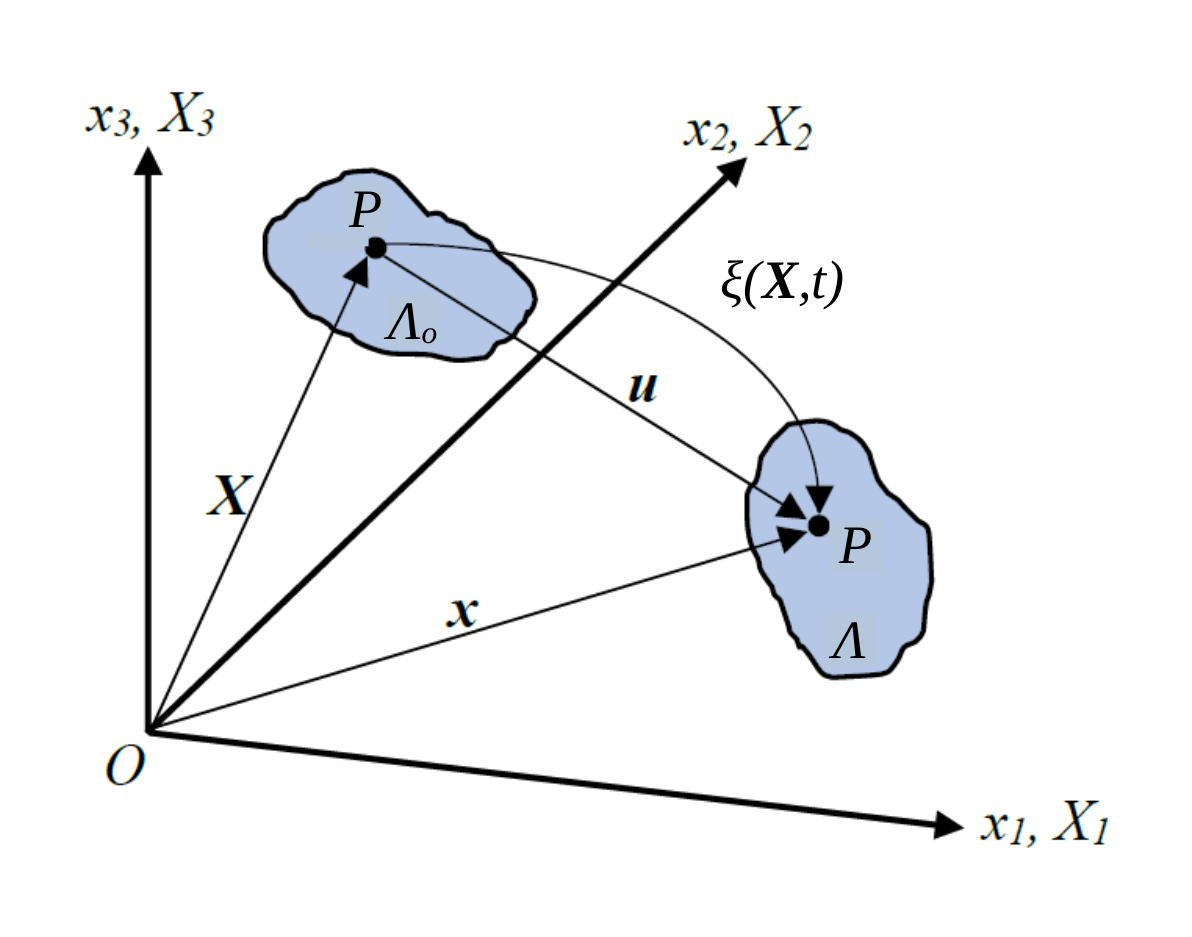}
	\caption{Description of the motion with respect to the initial and current configuration}
	\label{fig:motion_deformation}
\end{figure}

For clarity in our discussion, in this section we describe the motion and deformation of an arbitrary solid body~$\Omega \in \mathbb R^{3}$, as depicted in Figure ~\ref{fig:motion_deformation}. 
The current configuration of the body, denoted as $\Lambda$, is derived from the displacement and deformation relative to the initial configuration $\Lambda_o$. In our notation, we use uppercase $\vec X$ to represent vector coordinates in the initial configuration and lowercase $\vec x$ in the current configuration. The mapping between the initial and current configurations is denoted as~$\xi(\vec X, \, t)$. It is important to note that this mapping encompasses both the displacement and deformation between the two considered configurations. Furthermore, we specifically consider an infinitesimally small region~$P  \in \Omega$ within the domain with an initial position~$\vec{X}$ that occupies a new position~$\vec x$ in the current configuration. The relationship between the initial position~$\vec X$ and the current position~$\vec x$ is expressed as
\begin{equation} \label{eq:motion}
	\vec x = \vec X + \vec u \, ,
\end{equation}
where the vector $\vec{u}$ denotes the displacement of $P$. Utilising this relationship, we obtain the deformation gradient tensor $\vec F$ as 
\begin{equation} \label{eq:deformGrad}
	\vec{F} = \frac{\D \vec{x}}{\D \vec{X}} = \vec{I} + \frac{\D \vec{u}}{\D \vec{X}} \, .
\end{equation}
Based on the deformation gradient mentioned above, the stress tensor is calculated using a specific constitutive relation, which is dependent on the material type. In our examples, we consider a nearly incompressible hyperelastic Neo-Hookean constitutive model for rubber-like materials. As documented in~\cite{Lee2013}, the first Piola-Kirchhoff stress is decomposed into the summation of the volumetric and deviatoric components
\begin{equation}
    \vec{P} = \vec{P}_{vol} + \vec{P}_{dev} \, .
\end{equation}
Here, the volumetric and deviatoric stresses are calculated as follows
\begin{subequations}
	\begin{align}
    	\vec{P}_{vol} &= p J \vec{F}^{-\trans} \\
    	\vec{P}_{dev} &= \mu J^{-2/3} \left( \vec{F} - \frac{1}{3}(\vec{F} : \vec{F}) \vec{F}^{-\trans} \right ) \, .
    \end{align}
\end{subequations}
In these expressions, $p=\kappa(J-1)$ represents the hydrostatic pressure, $\mu$ is the shear modulus, $\kappa$ is the bulk modulus, and the Jacobian~$J$ is the determinant of the deformation gradient tensor, i.e., $J = \det(\vec F)$. The relation between the first Piola-Kirchhoff stress $\vec P$ and the Cauchy stress $\vec \sigma$ is established by 
\begin{equation}
	\label{eq:piola-kirchoff}
	\vec{\sigma} = \frac{1}{J} \vec{P} \cdot \vec{F}^{\trans} \, .
\end{equation}
 Note that the continuum assumption is employed throughout the derivation in this section. The discretisation of the continuous fields using particles, as described in Section~\ref{sec:discretisation}, will be discussed in the upcoming sections. 

\subsection{Total Lagrangian formulation of conservation laws}
\label{sec:conserv-laws}

In the total Lagrangian formulation, the initial configuration is used to describe the motion and deformation of the solid body. The conservation of linear momentum relative to the initial configuration is described as
\begin{equation}
	\label{eq:lin-momentum}
	\frac{\D \vec p}{\D t} -  \text{div}(\vec P) = \vec s_o \, .
\end{equation}
Here $\vec p$ denotes the momentum per unit of undeformed volume, 
$\vec s_o$ represents the source of force per unit of undeformed volume, 
and $\vec P$ is the first Piola-Kirchhoff stress. The subscript $o$ is used to indicate that the fields and their derivatives are associated with the reference, or undeformed, configuration $\Lambda_o$.  We discretise the stress field according to the SPH discretisation introduced in Section~\ref{sec:discretisation}, resulting in
\begin{equation}
	\label{eq:coptlf}
	 \frac{\D \vec p_i}{\D t} 
	= \vec s_{o, i} + \sum_{j \in \set N(\vec X_i)} V_{o, j} 
	  \vec P_j  {\nabla}_o \widetilde{W}_j^1 (\vec X_i) \, .
\end{equation}
Here $\vec X_j$ indicates the position of the $j-$th particle in the initial configuration, and $V_{o,j}$ denotes the undeformed volume of the $j-$th particle. In~\eqref{eq:coptlf}, we adopt a common approximation of derivative in SPH and the gradient of the first-degree corrected kernel $\nabla_o \widetilde{W}_j^1 (\vec X_i)$ is applied to ensure the zero-th and first order completeness. 

It is important to note at this stage that, following the correction, the value of the kernel gradient is no longer antisymmetric, i.e., ${\nabla}_o \widetilde{W}_{j}^1 (\vec X_i) \ne -{\nabla}_o \widetilde{W}_{i}^1 (\vec X_j)$. The direct implication of this property is that the interaction forces between two particles are not the same in the pairwise direction and, therefore, do not, in general, satisfy the conservation of linear momentum~\cite{frontiere2017}. In this work, we introduce $\vec{P}_i$ and modify the corrected kernel gradient to conserve the total linear momentum of the system, as also proposed in~\cite{frontiere2017}. The symmetrised conservation of linear momentum is expressed as follows
\begin{equation}
	\label{eq:symmetrised_com1}
	\frac{\D \vec{p}_i}{\D t} = \vec s_{o, i} + \frac{1}{2} \sum_{j \in \set N(\vec X_i)} V_{o,j} \left( \vec{P}_i + \vec{P}_j \right) \left(\nabla_o \widetilde{W}_j^1(X_i) - \nabla_o \widetilde{W}_i^1(X_j) \right)  \, .
\end{equation}
The formulations in~\eqref{eq:symmetrised_com1} ensure pairwise antisymmetric interaction forces, yielding linear momentum preservation. Details of the derivation can be found in ~\ref{sec:derivation_com}.

Following~\cite{BonetPart12015}, the additional relation used is the conservation of the deformation gradient $\vec F$, {which is derived from the geometrical compatibility.} The relation is governed by
{
\begin{equation}
	\frac{\D \vec F}{\D t} - \text{div}  \left( \frac{1}{\rho_o} \vec p \otimes \vec I \right) = \vec 0 \, .
	\label{eq:deform-grad-evo}
\end{equation} }
%
The discrete form of~\eqref{eq:deform-grad-evo} is again obtained using the SPH discretisation,
\begin{equation} 
	\label{eq:Fave}
	\frac{\D \vec F_i}{\D t} 
	= \sum_{j \in \set N(\vec X_i)} \frac{V_{o, j}}{\rho_{o, j}} 
	\vec p_j \otimes \nabla_o \widetilde{W}_{j}^1 (\vec X_i) \, .
\end{equation}
We reiterate here that the use of the gradient of the first order corrected kernel $\nabla_o \widetilde{W}_{j}^1 (\vec X_i)$ ensures that linear field can be reproduced and its respective gradient can be exactly computed up to machine precision~\cite{frontiere2017}. 


In our solid dynamics simulations, we use the Jameson-Schmidt-Turkel (JST) method~\cite{AguirreJST2014, Lee2016} to suppress the presence of excessive dissipation commonly observed in cases involving large strain. The JST method introduces an additional term~$\vec D(\vec p)$ in terms of the momentum field $\vec p$ into the discretised conservation of momentum, resulting in,
\begin{equation}
	\label{eq:cop}
	\frac{\D \vec p_i}{\D t} 
	= \vec s_{o, i} + \frac{1}{2} \sum_{j \in \set N(\vec X_i)} V_{o, j} \left(
	\vec P_i + \vec P_j \right) \left( \nabla_o \widetilde{W}_j^1 (\vec X_i)- 
	\nabla_o \widetilde{W}_i^1 (\vec X_j) \right) + \vec D(\vec p)\, ,
\end{equation}
where $\vec D(\vec p)$ consists of two terms involving harmonic and biharmonic operators
\begin{equation}
	\label{eq:lap_D}
	\vec D(\vec p) = 
	            \eta^{(2)} C_p h \sum_{j \in \set N(\vec X_i)} V_{o, j} 
	            (\vec p_j - \vec p_i) \widetilde{\nabla}_o^2 W_j 
	            (\vec X_i) - \eta^{(4)} C_p h^3 \sum_{j=1}^{N_i} 
	            V_{o, j} (\nabla_o^2 \vec p_j - \nabla_o^2 \vec p_i) 
	            \widetilde{\nabla}_o^2 W_j (\vec X_i) \, .        
\end{equation}
The constant $C_p$ denotes the speed of $p$-wave travelling through the material~\cite{Lee2022} and the parameter $h$ denotes the minimum smoothing length within the entire domain. In most of our dynamic simulations, we use the constants~$\eta^{(2)} = 0$ and $\eta^{(4)} = 0.125$ according to~\cite{Lee2016}. 

Further inspection of the equation~\eqref{eq:lap_D} indicates that the Laplacian term must satisfy completeness. For the second derivative of the kernel to reproduce a constant function, the kernel must satisfy second-order completeness, which is often deemed computationally prohibitive. In this work, we adopt the correction for the Laplacian term following the approach in~\cite{BonetKulasegaram}, resulting in the term~$\widetilde{\nabla}_o^2 W_j (\vec X_i)$. 
The resulting corrected Laplacian of the kernel sufficiently ensures that the second derivative is zero for linear and constant fields, and constant for a quadratic field. To evaluate the term $\nabla_o^2 \vec p_i$ in~\eqref{eq:lap_D}, we expand with SPH discretisation involving the corrected Laplacian~$\widetilde{\nabla}_o^2 W_j (\vec X_i)$, i.e., 
\begin{equation}
	\nabla_o^2 \vec p_i= \sum_{j \in \set N(\vec X_i)} V_{o, j} \left( \vec p_j - \vec p_i \right) \widetilde{\nabla}_o^2 W_j (\vec X_i) \, .
\end{equation}
Furthermore, we note that the corrected Laplacian of the kernel is not symmetric, i.e.,  $\widetilde{\nabla}_o^2 W_j (\vec X_i) \neq \widetilde{\nabla}_o^2 W_i (\vec X_j)$, which can affect the preservation of linear momentum. To address this concern, we use the symmetrised form for each particle pair, leading to the final expression of $\vec D(\vec p)$
\begin{equation}
	\begin{aligned}
	\label{eq:lap_D_sym}
	\vec D(\vec p) = \, 
	&\frac{1}{2} \eta^{(2)} C_p h \sum_{j \in \set N(\vec X_i)} V_{o, j} 
	\left( \vec p_j - \vec p_i \right) \left( \widetilde{\nabla}_o^2 W_j (\vec X_i) + \widetilde{\nabla}_o^2 W_i (\vec X_j)\right)  - \\
	&\frac{1}{2} \eta^{(4)} C_p h^3 \sum_{j \in \set N(\vec X_i)} V_{o, j} 
	\left( \nabla_o^2 \vec p_j - \nabla_o^2 \vec p_i \right) 
	\left( \widetilde{\nabla}_o^2 W_j (\vec X_i) + \widetilde{\nabla}_o^2 W_i (\vec X_j)\right) \, .
	\end{aligned}
\end{equation}
%
%

Our final component in formulating the total Lagrangian SPH involves time integration. To illustrate our strategy, we define the rate of change of a generic particle property $\vec U$ as
\begin{equation}
	\frac{\D \vec U}{\D t} = R(\vec U, \, t) \, .
\end{equation}
An example of the particle property $\vec U$ is the linear momentum $\vec p$. In this case, the above expression resembles equation~\eqref{eq:cop}, where here $R$ indicates the right hand side operator. In this work, we employ the explicit three-stage time integrator to update the particle property~$\vec U$ at every time step~\cite{Lee2022},  
\begin{subequations}
	\label{eq:time-int}
	\begin{align}
	\vec U^{*} &= \vec U^n + R(\vec U^n, \, t^n)  \, \Delta t \,  \\
	\vec U^{**} &= \frac{3}{4} \vec U^n + \frac{1}{4} \left( \vec U^* + R(\vec U^*, \, t^n)  \, \Delta t \right) \,  \\
	\vec U^{n+1} &= \frac{1}{3} \vec U^n + \frac{2}{3} \left( \vec U^{**} + R(\vec U^{**}, \, t^n)  \, \Delta t \right) \, .		
	\end{align}
\end{subequations}
The time increment $\Delta t = t^{n+1} - t^n$ is calculated based on the 
Courant-Friedrichs-Lewy (CFL) condition
\begin{equation}
	\label{eq:cfl}
	\Delta t = \alpha_{CFL} \cdot min \left(\frac{r_{ij}}{C_{p,ij}} \right) \, ,
\end{equation}
where $r_{ij}=\| \vec{x}_i - \vec{x}_j \|$ is the distance between particle $i$ and $j$ with respect to the current configuration, and $C_{p,ij}$ is the average $p$-wave speed on the interaction pair. The constant $\alpha_{CFL}$ is commonly case dependent, but our default setting is as high as $\alpha_{CFL}=0.9$. {According to the reference~\cite{Lee2022}, $C_{p,ij}$ is calculated as follows
\begin{equation}
	C_{p,ij} = \frac{C_p^{Lin}}{\lambda_{ij}^{Ave}};
	\hspace{1cm} 
	C_p^{Lin} = \sqrt{\frac{\lambda+2\mu}{\rho}};
	\hspace{1cm}
	\lambda_{ij}^{Ave}=\frac{1}{2} \left( \lambda_i + \lambda_j \right) \, ,
\end{equation}
where $\mu$ is the shear modulus, $\lambda$ is the Lam\'{e} first parameter, and $\lambda_i$ and $\lambda_j$ are the minimum stretch of particles $i$ and $j$, respectively.}

%% file: pseudostaggered.tex
\section{Alternative formulation for computing deformation gradient}
\label{sec:pseudo-staggered}


Spurious non-physical oscillations, namely the zero-energy mode or hourglass mode, have been observed in various computational methods including finite difference, finite elements, and particle method like SPH. These oscillations can be understood as the manifestation of displacement without any corresponding strain. In the context of SPH, this phenomenon arises due to the inaccurate  computation of gradients of the field, allowing for zero gradients not strictly associated with constant fields. Motivated by this insight, we elaborate on our proposed alternative for computing the deformation gradient with improved accuracy to suppress the presence of the  zero-energy mode. In this section, we begin by outlining our approach in a one-dimensional setting and subsequently integrate this idea into the total Lagrangian SPH framework in higher dimensions. 

\subsection{One-dimensional illustration}
For clarity and simplicity of the discussion, we illustrate our approach in a one-dimensional setting. Consider the domain $\Omega \in \mathbb R^1$ in the reference configuration, where we assume that the deformation field can be obtained at a finite set of particles with coordinates $X_j \in \Omega$. The deformation gradient at a specific particle $X_i$ can be obtained using the SPH discretisation scheme as discussed in Section~\ref{sec:discretisation},
\begin{equation}
	\label{eq:1d-df}
	\frac{\D x}{\D X} (X_i) \approx \sum_{j \in \set N(\vec X_i)}  \frac{\D W(X_i - X_j , h)}{\D x} \, V_j  \, x_j \, .
\end{equation}
The notation $\set N(X_i)$ denotes the list of particle indices belonging to the neighbourhood of $X_i$. The observation of zero-energy, or hourglass, mode in SPH is attributed to the discrete approximation of the deformation gradient, which permits it to be zero not strictly for a constant deformation field. This observation motivates the development of methods that redefine the gradient computation in SPH using an additional set of particles, known as stress points~\cite{Dyka1997, Vignjevic2000, Randles2000SP}. 

\begin{figure}[]
	\centering
	\includegraphics[width=1\textwidth]{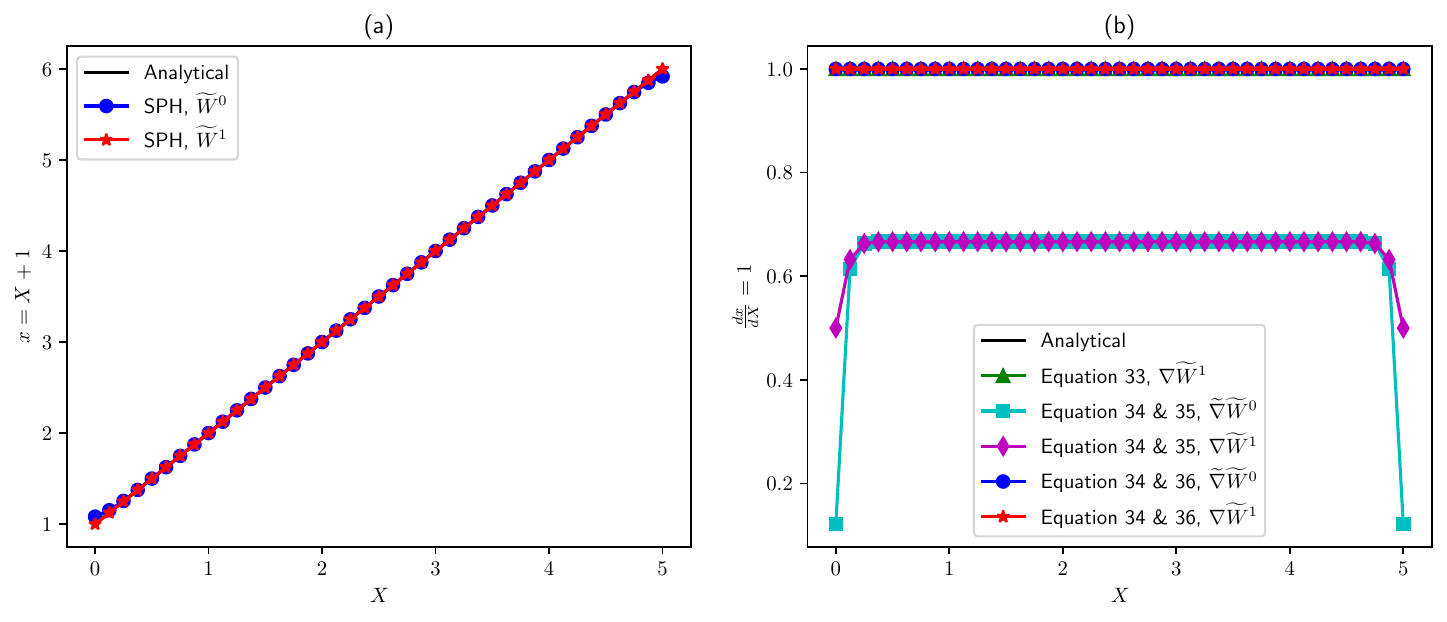}
	\caption{Field reconstruction from point values: (a) deformation , (b) deformation gradient.}
	\label{fig:grad-linear}
\end{figure}

Another possible formulation of the field derivative can be defined using the gradient of deformation obtained at the particles and the kernel $W$, in contrast to the standard SPH formulation in~\eqref{eq:1d-df}, i.e., 
\begin{equation}
	\label{eq:alternative-df}
	\frac{\D x}{\D X} (X_i)  \approx \sum_{j \in \set N(\vec X_i)} V_j W(X_i - X_j, \, h) \left( \frac{\D x}{\D X}  \right)_{ij} \, .
\end{equation}
This alternative formulation can be traced back to the application of a differential operator $D$ to the convolution between two continuous functions $W$ and $x$, which leads to $D\, (W \ast x) = D \, W \ast x = W \ast D \, x$. Note that in our one dimensional illustration, the operator $D$ is $\D / \D X$. A geometric interpretation of \eqref{eq:alternative-df} is that the deformation gradient at particle $X_i$ is approximated by a weighted average of the gradient measured between $X_i$ and its neighbour $X_j$. This pairwise gradient can be simply defined as
\begin{equation}
	\label{eq:alternative-df2}
	\left( \frac{\D x}{\D X}  \right)_{ij} = \frac{x_j - x_i}{X_j - X_i} \, .
\end{equation}
The above equation provides an estimate for the derivative between $X_i$ and the neighbouring particles at $X_j$, implicitly requiring the index $j \neq i$ to avoid a zero denominator. This approach can also be viewed as a \textit{bond-based} gradient estimation, as proposed in peridynamics~\cite{Tupek2014, Silling2017, Breitzman2018}. While the presented formulation in~\eqref{eq:alternative-df2} is straightforward, it lacks precision when approximating the gradient, even for a linear field $x = X + 1$, as illustrated in Figure~\ref{fig:grad-linear}. Notably, even in regions distant from the boundary, the estimated gradient does not equate to 1, stemming from the omission of the gradient estimate at $X_i$. 

Building on the aforementioned observation, we introduce an improved formulation of the bond-based gradient $\left(\dfrac{\D x}{\D X}\right)_{ij}$ by considering the gradient estimate at particle $X_i$,
\begin{equation}
	\label{eq:modified-bond}
	\left( \frac{\D x}{\D X}  \right)_{ij} = 
		\begin{dcases}
		\frac{x_j - x_i}{X_j - X_i}	&  \text{if } \,  j \neq i \\
		\sum_{j \in \set N(\vec X_i)} \frac{\D W(X_i - X_j, \, h)}{\D X} V_j  x_j   & \text{if } \,  j = i \, 
	\end{dcases} \, . 
\end{equation}
Note that for $i = j$, we adopt the standard gradient estimate in SPH~\eqref{eq:1d-df}. This formulation can also be viewed as a correction of~\eqref{eq:1d-df} evaluated at $X_i$ by considering the pairwise deformation gradient with its neighbours. As a preliminary assessment, we consider a uniform particle distribution in one-dimensional domain $\Omega = (0, 5)$ with spacing $\Delta x = 1/8$. Our objective is to approximate the gradient of the deformation field $x = X + 1$ using the standard formula~\eqref{eq:1d-df}, the bond-based estimate~\eqref{eq:alternative-df2}, and the improved formulation~\eqref{eq:modified-bond}. Here we consider a cubic b-spline kernel with first order correction with half-width $R = 4 \, \Delta x$. As depicted in Figure~\ref{fig:grad-linear}(b), the bond-based estimate~\eqref{eq:alternative-df2} exhibits discrepancy with the expected value. However, the proposed improvement to the bond-based formulation, as described in~\eqref{eq:modified-bond}, yields a correct deformation gradient of 1 as also achieved by the standard SPH described in~\eqref{eq:1d-df}. 

\begin{figure}[ht!]
	\centering
	\includegraphics[width=1\textwidth]{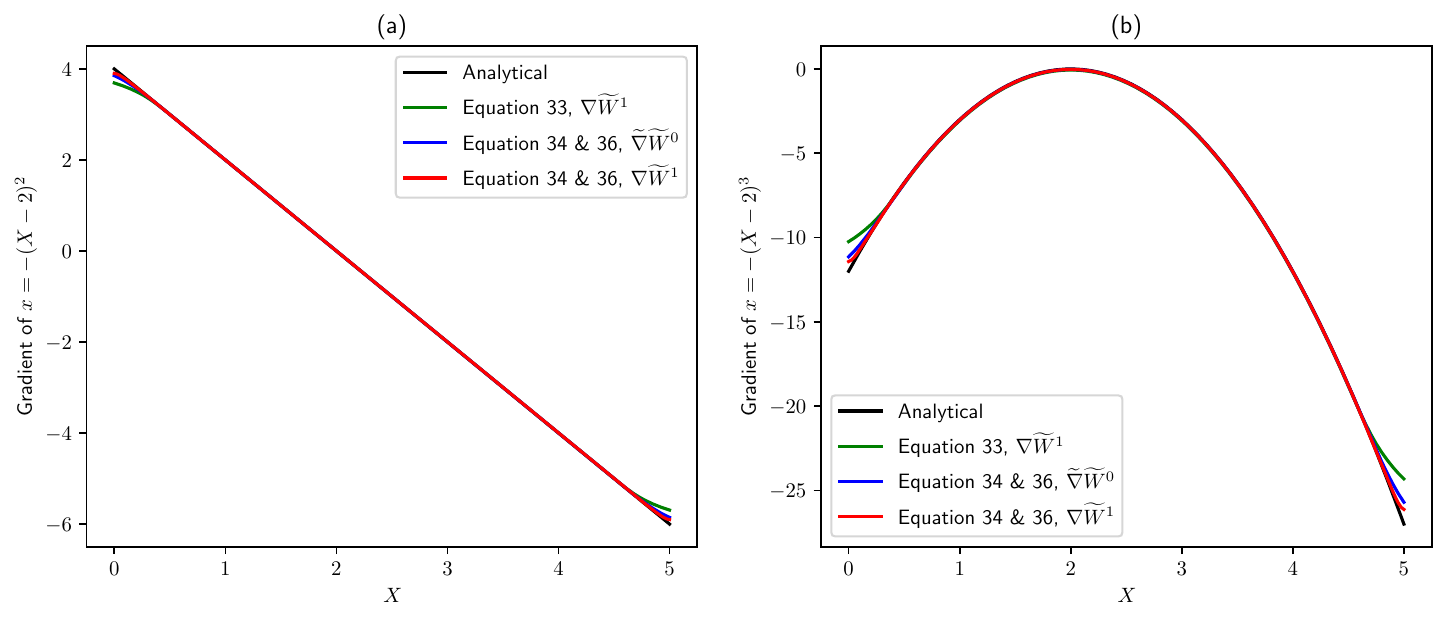}
	\caption{Gradient approximation of a deformation field from point values: (a) quadratic function, (b) cubic function}
	\label{fig:grad-field}
\end{figure}

A more detailed analysis is conducted by calculating the deformation gradient of higher-order deformation fields. Specifically, our focus lies on comparing the standard SPH formula~\eqref{eq:1d-df} using the gradient of first-order corrected kernel $\nabla \widetilde{W}^1$, and our proposed formula~\eqref{eq:modified-bond} with two variants of kernel gradients $\nabla \widetilde{W}^1$ and $\widetilde{\nabla} \widetilde{W}^0$. Here we use the notation $\nabla$ for a comparable referencing with Section~\ref{sec:corrected-sph}. The deformation gradients of quadratic and cubic polynomial deformation fields are presented in Figure ~\ref{fig:grad-field}. The results indicate that for higher polynomial orders, the standard SPH~\eqref{eq:1d-df} exhibits errors due to the maximum polynomial reproducibility being only up to linear. The elimination of this error requires satisfying higher-order conditions, demanding more intensive computational efforts. It is also evident from~\eqref{eq:modified-bond} that the proposed formulation significantly improves accuracy, especially in proximity to the boundary. 

\begin{figure}[ht!]
	\centering
	\includegraphics[width=1\textwidth]{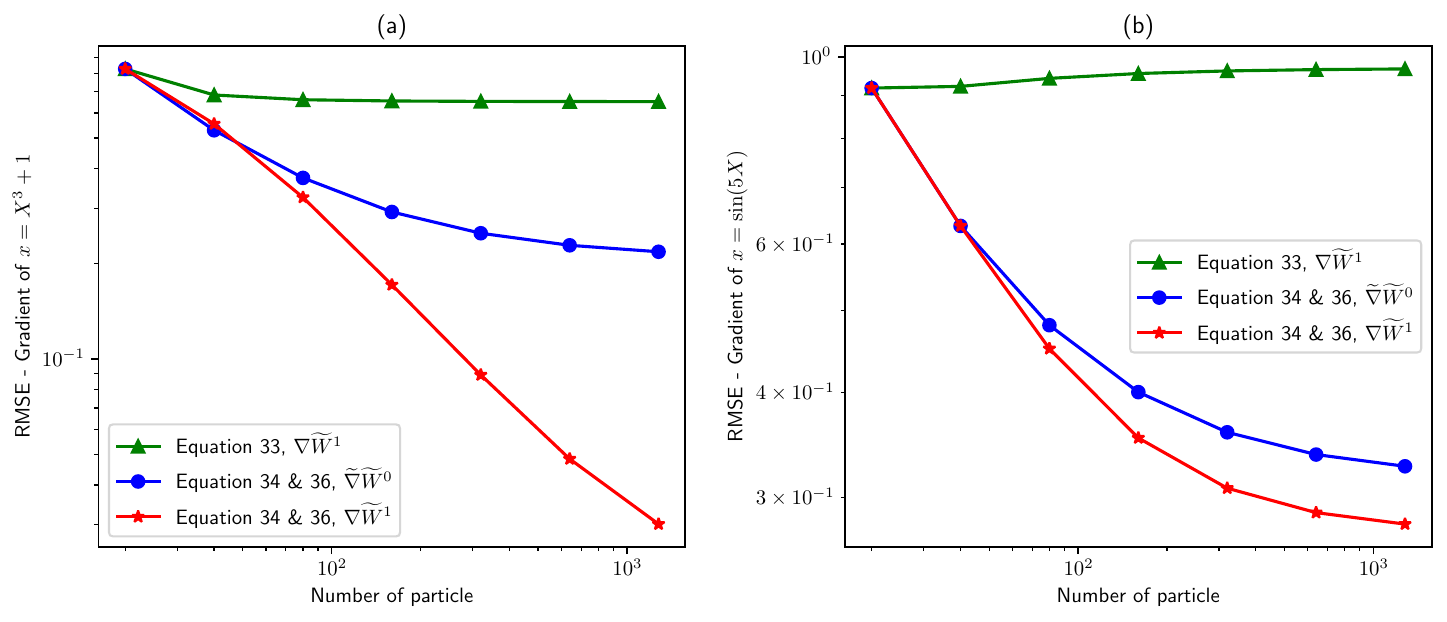}
	\caption{RMSE analysis of the SPH approximation: (a) gradient of cubic function, (b) gradient of sinusoidal function}
	\label{fig:rmse}
\end{figure}

Next, we assess the convergence of error in estimating the deformation gradient with respect to spatial refinement. As the deformation field, we employ cubic $x = X^3 + 1$ and periodic $x = \sin( 5 X)$ functions. Figure~\ref{fig:rmse} shows the root mean square error (RMSE) of the deformation gradient with respect to the analytical solutions. {As the number of particles increases, we observe that the RMSE of the gradient approximation of a cubic deformation field decreases. In contrast, when approximating the gradient of a sine deformation field, the RMSE plateaus at a certain threshold. This can be explained by the fact that the Taylor expansion of a sine function results in an infinite series of polynomials of ascending order. The accumulation of higher-order approximation errors leads to a higher error threshold in the sine case, while the error threshold remains much lower for the cubic case.} Based on the results obtained, it is evident that the RMSE of classical SPH~\ref{eq:1d-df} does not decrease as $h$ is refined and more particles are involved. In contrast, reduction of error is observed for the proposed formulation~\ref{eq:modified-bond}, with the one using $\nabla \widetilde{W}^1$ exhibiting lower error than the one using $\widetilde{\nabla} \widetilde{W}^0$. This discrepancy is less emphasised when dealing with oscillatory functions. This justifies that the proposed \text{improved bond-based} formulation significantly improves the estimation of the deformation gradient from the field values obtained at the particles. 

\subsection{SPH implementation of the proposed deformation gradient}

In this section we extend the one-dimensional formulation of the improved bond-based deformation gradient~\eqref{eq:modified-bond} to higher dimensions and describe its implementation in our TLSPH framework, as described in Section~\ref{sec:tlsph}. As a preliminary reminder, in our TLSPH, the deformation gradient at each particle $\vec F_i$ is computed at every time step by first evaluating the conservation of the deformation gradient~\eqref{eq:Fave}, and subsequently integrating over a time step $\Delta t$. Following the evaluation of $\vec F_i$, the modified deformation gradient $\vec F_i^c$ is then obtained analogous to the one-dimensional expression~\eqref{eq:alternative-df},
\begin{equation}
	\label{eq:nd-df}
	\vec F_i^c = \sum_{j \in \set N(\vec X_i)} V_j W(\vec X_i - \vec X_j, h) \vec F_{ij} \, .
\end{equation}

One possible expression for the bond-based deformation gradient $\vec F_{ij}$ can be obtained purely from kinematics information, serving as a direct extension from our one-dimensional illustration, 
\begin{equation}
	\label{eq:modified-bond-nd}
	\vec F_{ij} = 
	\begin{dcases}
		\frac{\vec x_j - \vec x_i}{r_{o,ij}}	\otimes \frac{\vec X_j - \vec X_i}{r_{o,ij}}	 &  \text{if } \,  j \neq i \\
		\vec F_i  & \text{if } \,  j = i
	\end{dcases} \, .
\end{equation}
The term $r_{o, ij}$ indicates the Euclidian distance between particle $i$ and $j$ in the initial configuration, i.e., $\| \vec X_j - \vec X_i \|$. Note that this expression is still consistent with~\eqref{eq:modified-bond} since $\frac{\vec X_j - \vec X_i} {r_{o,ij}} = 1$ in one-dimension. The inter-particle component $\frac{\vec x_j - \vec x_i} {r_{o,ij}} \otimes \frac{\vec X_j - \vec X_i} {r_{o,ij}}$ is rank deficient. An alternative expression for the inter-particle deformation gradient (i.e., when $j \neq i$) is proposed in~\cite{Breitzman2018} by correcting the \textit{averaged} deformation gradient $\vec F_i$ projected to each bond with the actual deformation~\eqref{eq:modified-bond-nd}, which reads
\begin{equation}
	\label{eq:bond-pred-corr}
	\vec F_{ij} = 
	\begin{dcases}
	\vec F_i  - \left( \left( \vec F_i \cdot \frac{\vec X_j - \vec X_i} {r_{o,ij}} \right)  \otimes \frac{\vec X_j - \vec X_i} {r_{o,ij}} \right) +
	\left( \frac{\vec x_j - \vec x_i} {r_{o,ij}} \otimes 
	\frac{\vec X_j - \vec X_i} {r_{o,ij}} \right)	 &  \text{if } \,  j \neq i \\
		\vec F_i  & \text{if } \,  j = i
	\end{dcases} \, . 
\end{equation}

{In~\eqref{eq:bond-pred-corr}, the individual bond-based deformation gradient between the pair $i - j$ is obtained by first correcting the averaged deformation gradient $\vec F_i$ with the component derived from mapping a material line describing a bond between the pair $i - j$, i.e., $\frac{\vec X_j - \vec X_i}{r_{o,ij}}$, with $\vec F_i$. The actual deformation gradient $\left( \frac{\vec x_j - \vec x_i} {r_{o,ij}} \otimes \frac{\vec X_j - \vec X_i} {r_{o,ij}} \right)$ is then added to obtain the actual deformation gradient in each bond. As an illustrative example, we consider a case with unidirectional deformation, specifically in the $x$-direction for a bond between particle $i-j$ both aligned along the $x$-axis. The application of~\eqref{eq:modified-bond-nd} in this scenario yields  
\begin{equation}
	\label{eq:fij_full}
	\vec F_{ij} = 
	\left[ {\begin{array}{cc}
	\frac{(x_j-x_i)(X_j-X_i)}{r_{o,ij}^2} & 0. \\
	0. & 0. \\
	\end{array} } \right] \, ,
\end{equation}
}
{where $\vec F_{ij}$ is evidently rank deficient. If we now assume a deformation gradient $\vec F_i = \left[ {\begin{array}{cc}
	F_{i,xx} & 0. \\
	0. & 1. \\
	\end{array} } \right]$, then using~\eqref{eq:bond-pred-corr} we obtain a full-rank deformation gradient,
\begin{equation}
	\label{eq:fij_full}
	\vec F_{ij} = 
	\left[ {\begin{array}{cc}
	F_{ij,xx} & 
	0. \\
	0. & 
	1. \\
	\end{array} } \right] \, ,
\end{equation}
where
$F_{ij,xx} = F_{i,xx} - \frac{F_{i,xx}(X_j-X_i)(X_j-X_i)}{r_{o,ij}^2} + \frac{(x_j-x_i)(X_j-X_i)}{r_{o,ij}^2}$, and $F_{i,xx} = \left( 1+\frac{\partial u_x}{\partial x} \right)$.} 
Our final note concerns the computation of the particle deformation $\vec x_i$. In this work we obtain $\vec x_i$ from the momentum update~\eqref{eq:cop} and subsequently performing time integrations.

%% file: examples.tex
\section{Numerical examples}

We present in this section several numerical experiments with increasing complexity to investigate the convergence property and the performance of the improved bond-based deformation gradient formulation proposed in this work. Our first analysis focuses on the convergence of $L_2$-norm and $H_1$-seminorm errors in the displacement field using a two-dimensional swinging plate example. Within this study, we also investigate the effect of the kernel width as well as the correction type on the convergence. Additionally, we study the preservation of energy, linear momentum, and angular momentum through various two- and three-dimensional problems. These examples demonstrate that the proposed method improves numerical stability, enabling computations involving larger strains. 

\subsection{Convergence analysis using swinging plate problem}
\label{sec:swinging-plate}
\begin{figure}[]
	\centering
	\includegraphics[width=0.4\textwidth]{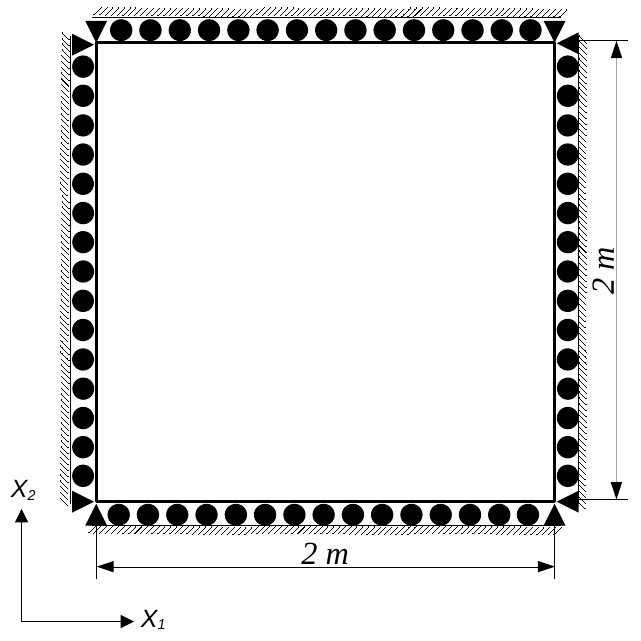}
	\caption{Configuration and dimension of the swinging plate for the convergence test}
	\label{fig:swinging}
\end{figure}

First, we investigate the convergence of the solution to the swinging plate problem~\cite{Lee2016,Lee2017,Lee2019,Park2019} in $L_2$-norm and $H_1$-seminorm errors of the displacement field. The domain is a square plate $\Omega \in \mathbb{R}^2$ with a width of 2~m. We consider a rubber-like material characterised by a nearly-incompressible Neo-Hookean constitutive model with a uniform Young's modulus $E = 17$~$\text{MPa}$, Poisson's ratio $\nu = 0.495$, and density $\rho_o = 1100$ $\text{kg/m}^3$. The displacement normal to the boundaries is constrained, as illustrated in Figure~\ref{fig:swinging}. An initial velocity is applied to the system at time $t = 0$, which is derived from a prescribed solution
\begin{equation}
	\vec{u}(\vec{x},t) = 
	U \sin{(\omega t)}
	\left \{ {\begin{array}{c}
			- \sin{(\frac{\pi}{2}X_1)} \cos{(\frac{\pi}{2}X_2)}  \\
			\cos{(\frac{\pi}{2}X_1)} \sin{(\frac{\pi}{2}X_2)}  \\
	\end{array} } \right \}  \, ,
	\label{eq:swinging_displacement}
\end{equation}
where $\omega$ represents the angular frequency given by
\begin{equation}
	\omega = \frac{\pi}{2} \sqrt{2\mu / \rho_o} \, .
\end{equation}
The variable $\mu = \frac{E}{2(1+\nu)}$ is the shear modulus, and the amplitude is chosen as $U=0.01~\text{m/s}$. Analytical solutions for the velocity and displacement gradient can be obtained through differentiation with respect to both time and space. 

The $L_2$-norm error of the displacement is defined as follows,
\begin{equation}
	\| \vec{u} - \vec{u}_{h} \| = \left( \int (\vec{u} - \vec{u}_{h}) \cdot (\vec{u} - \vec{u}_{h}) \,  \D \Omega \right)^{\frac{1}{2}} \, ,
	\label{eq:l2norm}
\end{equation}
and the $H_1$-seminorm error is defined as
\begin{equation}
	| \vec{u} - \vec{u}_{h} | =  \left( \int \left(\frac{\partial \vec{u}}{\partial \vec{X}} - \frac{\partial \vec{u}_{h}}{\partial \vec{X}}\right) : \left(\frac{\partial \vec{u}}{\partial \vec{X}} - \frac{\partial \vec{u}_{h}}{\partial \vec{X}} \right) \, \D \Omega  \right)^{\frac{1}{2}} \, ,
	\label{eq:h1seminorm}
\end{equation}
where $\vec{u}_{h}$ represents the numerical displacements. The integrations in~\ref{eq:l2norm} and~\ref{eq:h1seminorm} are evaluated numerically using Gauss quadratures. Specifically, the plate is first subdivided into $10 \times 10$ tiles, and then $5 \times 5$ Gauss points are distributed within each tile.

\begin{figure}[ht!]
	\centering
	\includegraphics[width=0.9\textwidth]{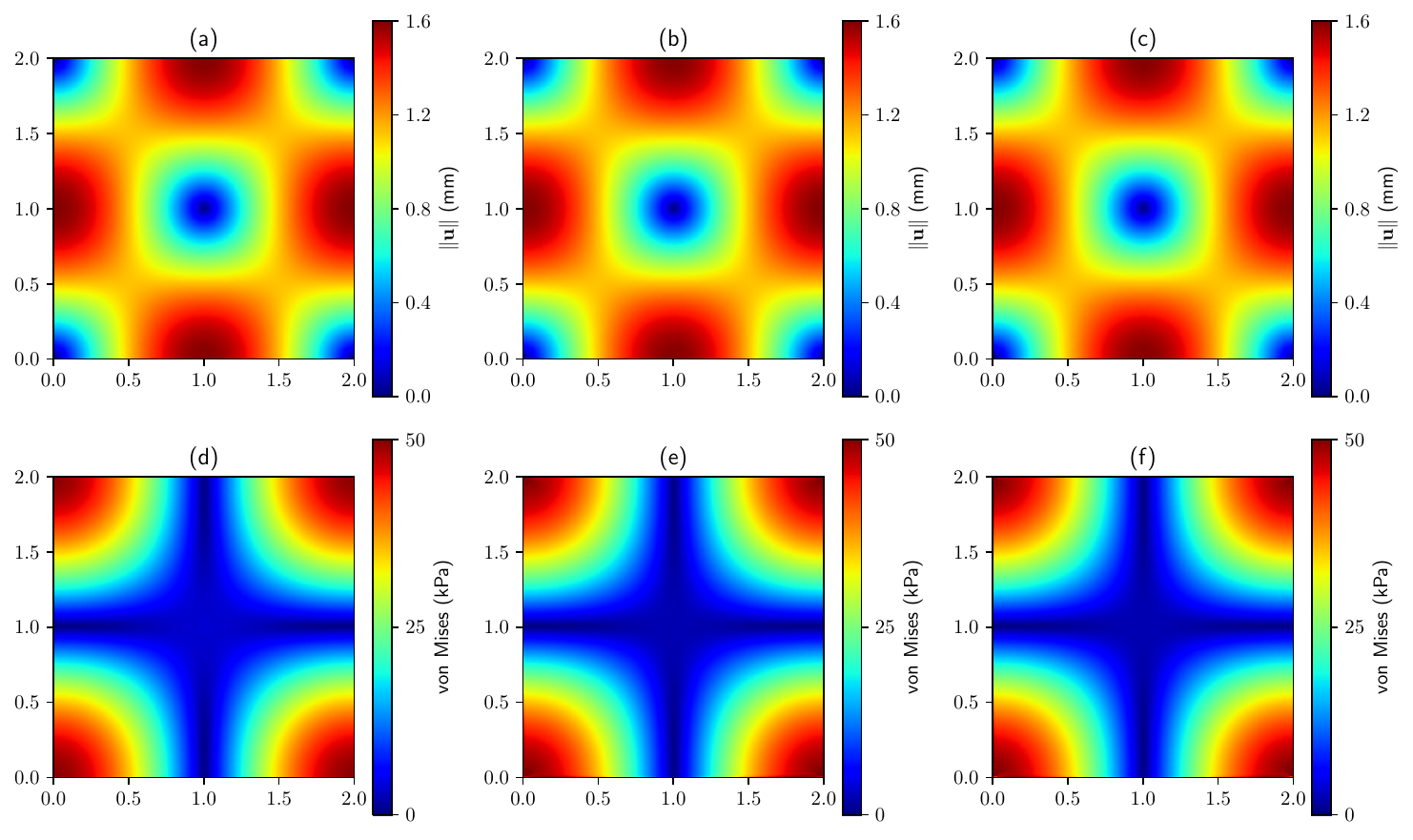}
	\caption{The displacement (m) and von Mises stress (Pa) of the swinging plate: (a) analytical displacement, (b) displacement using $\widetilde{W}^0_j(\vec{X}_i)$ \& $\widetilde{\vec{\nabla}}_0\widetilde{W}_j^0(\vec{X}_i)$, (c) displacement using $\widetilde{W}_j^1(\vec{X}_i)$ \& $\vec{\nabla}_0\widetilde{W}_j^1(\vec{X}_i)$, (d) analytical von Mises stress, (e) von Mises stress using $\widetilde{W}_j^0(\vec{X}_i)$ \& $\widetilde{\vec{\nabla}}_0\widetilde{W}_j^0(\vec{X}_i)$, (f) von Mises stress using $\widetilde{W}_j^1(\vec{X}_i)$ \& $\vec{\nabla}_0\widetilde{W}_j^1(\vec{X}_i)$}
	\label{fig:swinging_res}
\end{figure}

\begin{figure}[ht!]
	\centering
	\includegraphics[width=1\textwidth]{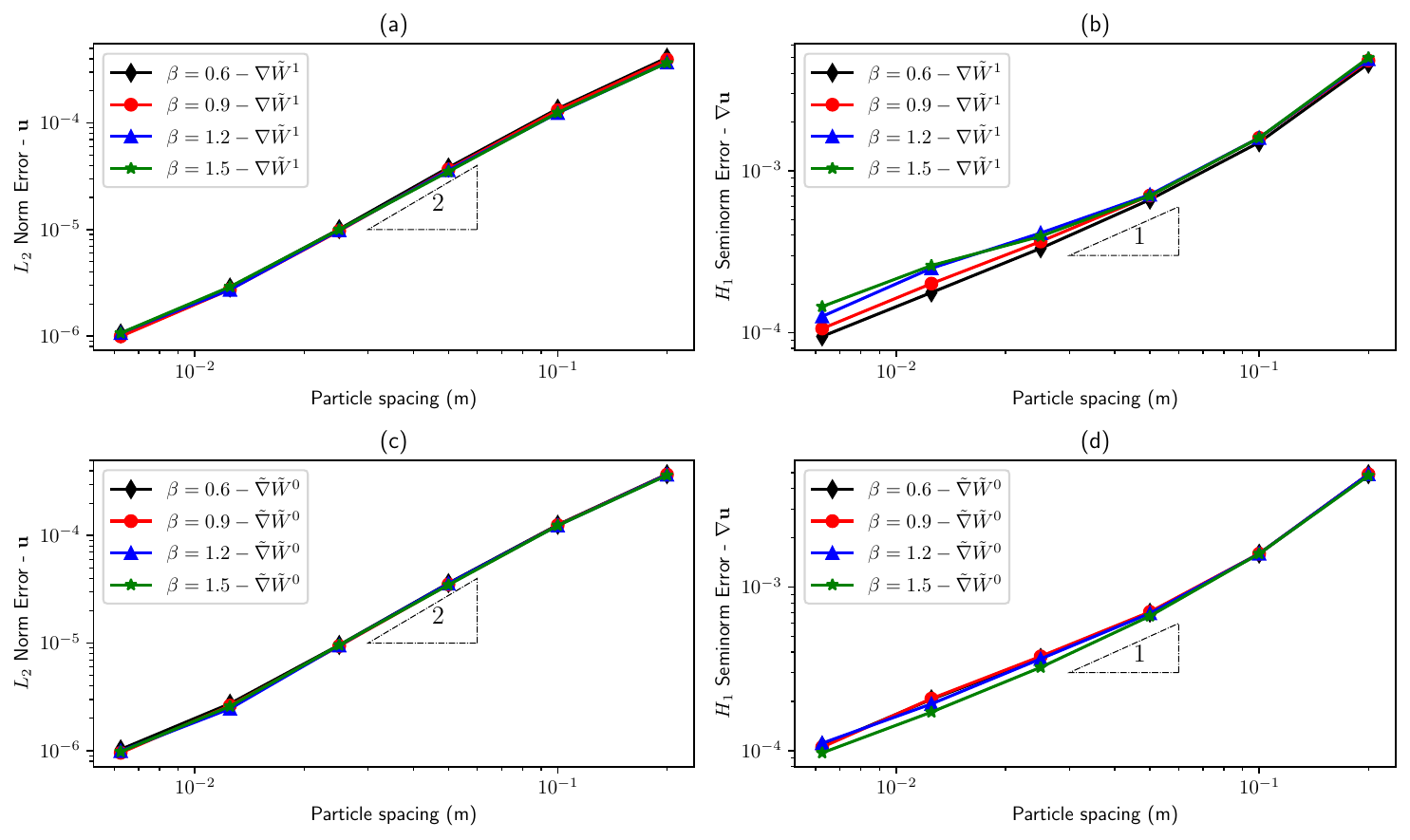}
	\caption{Convergence rate of the proposed SPH: (a) displacement $L_2$ norm error using $\widetilde{W}_i^1(\vec{X}_j)$ \& $\vec{\nabla}_0\widetilde{W}_i^1(\vec{X}_j)$, (b) displacement gradient $H_1$ seminorm error using $\widetilde{W}_i^1(\vec{X}_j)$ \& $\vec{\nabla}_0\widetilde{W}_i^1(\vec{X}_j)$, (c) displacement $L_2$ norm error using $\widetilde{W}^0_i(\vec{X}_j)$ \& $\widetilde{\vec{\nabla}}_0\widetilde{W}_i^0(\vec{X}_j)$, (d) displacement gradient $H_1$ seminorm error using $\widetilde{W}^0_i(\vec{X}_j)$ \& $\widetilde{\vec{\nabla}}_0\widetilde{W}_i^0(\vec{X}_j)$}
	\label{fig:swing_l2h1}
\end{figure}

In this example, we also study the effect of kernel correction on the obtained numerical solutions. Specifically, we consider the zero-th order corrected kernel~$\widetilde{W}^0_j(\vec{X}_i)$ and its corrected gradient~$\widetilde{\vec{\nabla}}_0 \widetilde{W}^0_j(\vec{X}_i)$, as well as the first order corrected kernel~$\widetilde{W}^1_j(\vec{X}_i)$ and its gradient~$\vec{\nabla}_0 \widetilde{W}^1_j(\vec{X}_i)$, as discussed in Section~\ref{sec:corrected-sph}. We use a cubic spline kernel with size defined by a parameter $\beta=h/d$, representing the ratio between the smoothing length $h$ and the particle spacing $d$. In our simulations, we vary $\beta$ within the set $\{ 0.6,0.9,1.2,1.5 \}$. The displacement and stress contours are shown in Figure~\ref{fig:swinging_res}. Qualitatively, it is evident that the contours are smooth without non-physical oscillations. Furthermore, the numerical results agrees well with the analytical solutions presented in tiles (a) and (d). The $L_2$-norm and $H_1$-seminorm errors are shown in Figure ~\ref{fig:swing_l2h1}. Here we consider an initial particle configuration with uniform spacing $d = 0.2$ and successively refine the spacing by a factor of 0.5. The two JST parameters in this example are $\eta^{(2)}=0$ and $\eta^{(4)}=0.125$. It is evident that the proposed method yields the optimal convergence rate of 2 in $L_2$-norm and 1 in $H_1$-seminorm. In the $L_2$-norm, the two kernel-gradient pairs lead to the optimal convergence rate of 2, with minimal variations with respect to the parameter $\beta$. In the $H_1$-seminorm these variations are more emphasised but still roughly yield to convergence rates between 0.9 to 1.2.

%

\subsection{Momentum preservation analysis using spinning plate and cube problems}
\label{sec:spinning-plate}
\begin{figure}[]
	\centering
	\includegraphics[width=0.4\textwidth]{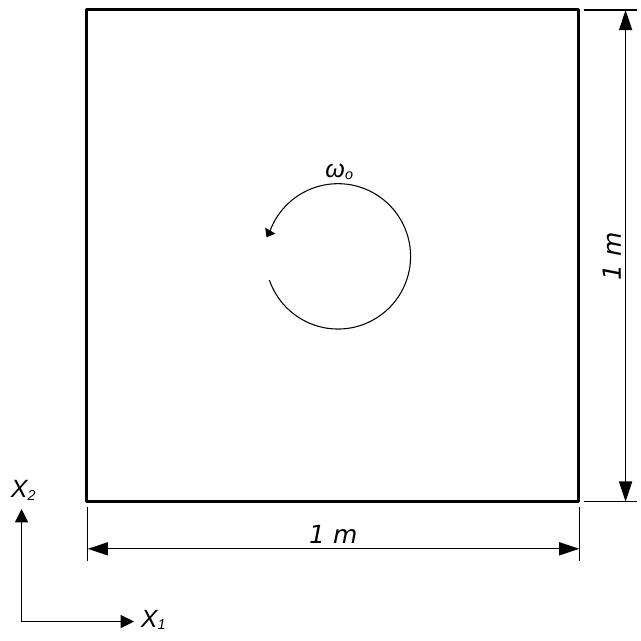}
	\caption{Configuration and dimension of the spinning plate}
	\label{fig:spinning_plate_specs}
\end{figure}

Next we investigate the preservation of linear and angular momentum using a spinning plate problem, following the example demonstrated in~\cite{Lee2014}, and aim to identify any spurious instability modes in the long-term response. We consider a square plate with width $1~\text{m}$ embedded in $\mathbb{R}^3$ with initial angular velocity vector $\vec{\omega}_o=(0,0,\omega_3)^\trans$ as depicted in Figure ~\ref{fig:spinning_plate_specs}. The axis of rotation is placed at the centre of the plate, with $\omega_3=105~\text{rad/s}$. The initial velocity can be computed as
\begin{equation}
		\vec{v}_o(\vec{X})=\vec{\omega}_o \times \vec{X} \,  ,
\end{equation}
where $\vec{X} = (X_1,X_2,0)^\trans$ represents the position in the reference configuration measured from the centre of rotation. We consider a rubber material characterised by a modulus of elasticity $E=17~\text{MPa}$, Poisson's ratio $\nu=0.45$, and density $\rho_o=1100~\text{kg/m}^3$. In this example, we discretise the plate using $200 \times 200$ uniformly distributed particles. We choose a cubic b-spline kernel with a smoothing length-to-particle spacing ratio  $\beta = 0.9$. Furthermore, in this example we employ the first order kernel correction and its gradient to ensure linear reproducibility.

\begin{figure}[]
	\centering
	\includegraphics[width=1\textwidth]{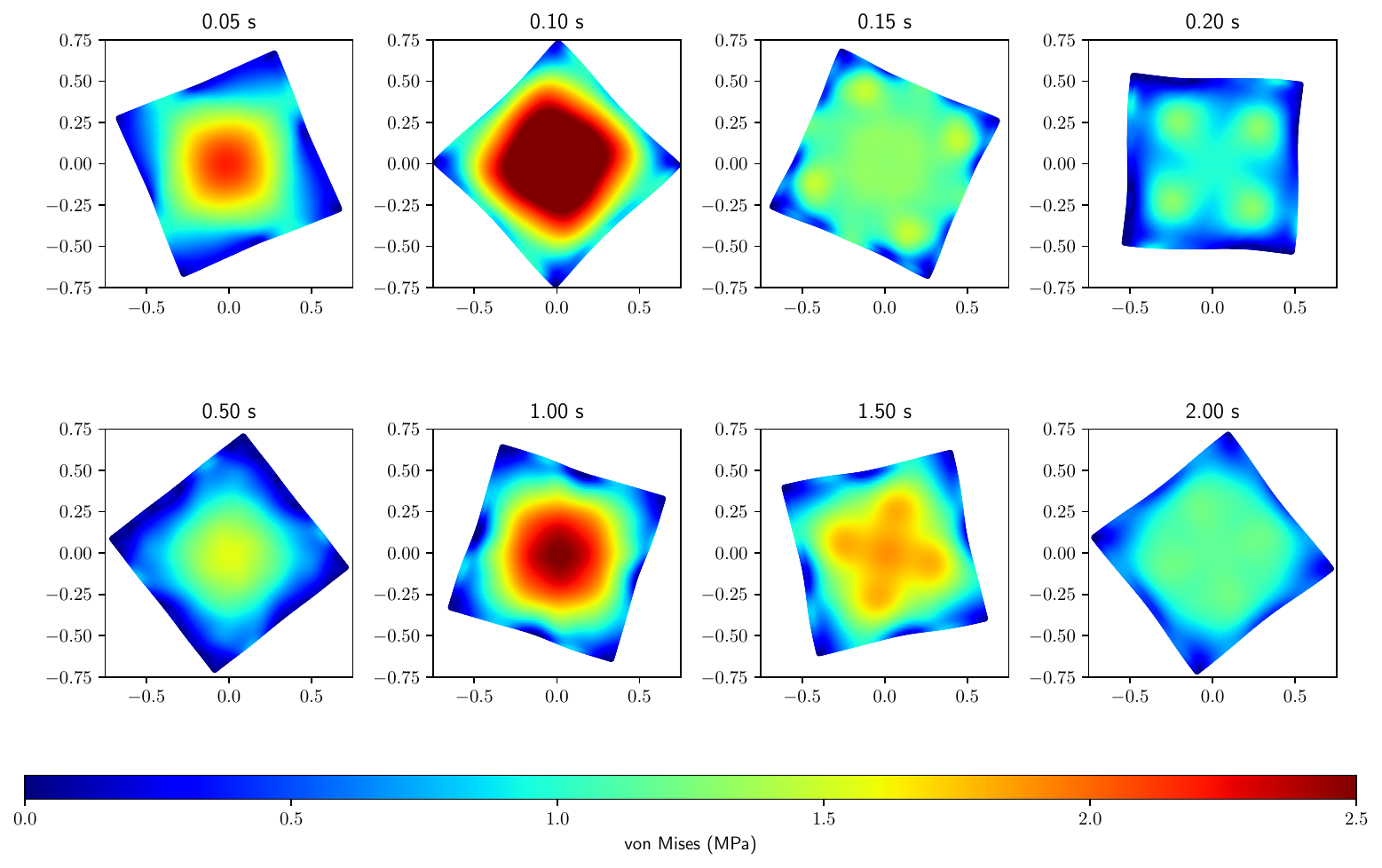}
	\caption{Stress contours of the spinning plate at various time frames}
	\label{fig:spinning_contour}
\end{figure}

We qualitatively examine the presence of the spurious oscillations by inspecting the von Mises stress distribution, as presented in Figure~\ref{fig:spinning_contour}. In the short term response, up to $t = 0.2~\text{s}$, depicted in the first row of Figure~\ref{fig:spinning_contour}, the deformation shape obtained through the proposed method aligns with the results documented in~\cite{Lee2014}. To further assess the performance of our proposed method, we extend the simulation to examine its long term response, reaching up to $t = 2~\text{s}$. As evident from the second row of Figure~\ref{fig:spinning_contour}, our proposed method endures the extended simulation duration, remaining free from spurious artifacts. This demonstrates stability and robustness over an extended period, reinforcing the efficacy of the proposed approach.

\begin{figure}[]
	\centering
	\includegraphics[width=0.8\textwidth]{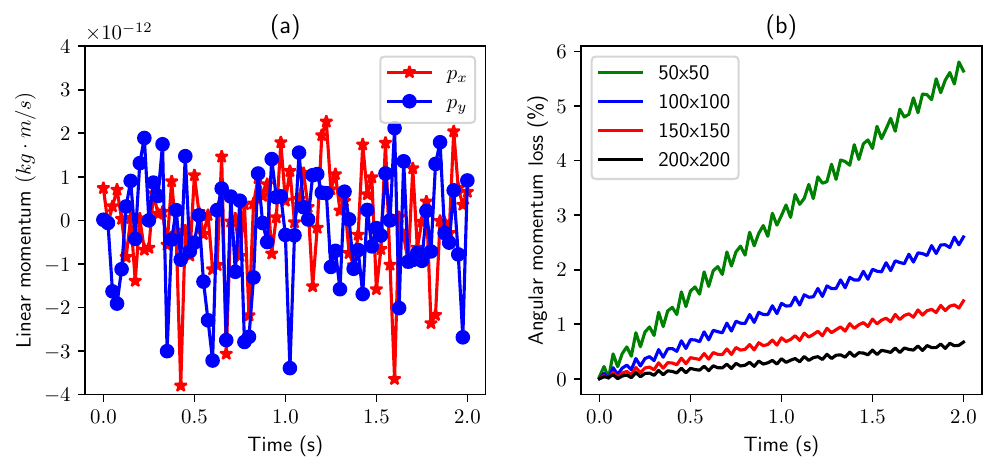}
	\caption{Linear and angular momentum of the spinning plate problem: (a) linear momentum using $50 \times 50$ particles, (b) effect of particle refinement on angular momentum loss}
	\label{fig:spinning_momentum}
\end{figure}

Next we examine the total linear and angular momentum in the system. Given the pure rotational movement of the plate, the total linear momentum should be zero, and the total angular momentum should be determined by the total of moment of inertia of each particle in the plate multiplied by its angular velocity. These values should remain constant throughout the simulation. The total linear momentum for the spinning plate problem is presented in Figure~\ref{fig:spinning_momentum}(a). It is evident that our method ensures the preservation of linear momentum throughout the simulation time. 
As discussed in~\cite{BonetVariational1999}, the preservation of linear momentum is inherent to the kernel correction and its gradient being used. However, preserving angular momentum remains a challenge due to the corrected kernel being non co-linear, which introduces additional moments. {In this work, we analyse how angular momentum loss changes with respect to spatial refinement. As shown in  Figure~\ref{fig:spinning_momentum}(b), the coarsest particle arrangement yields the highest loss of around $6\%$ at $t = 2$~s, whereas the lowest error is observed when the particle density is highest. This effect occurs due to shorter moment arms between pairs of particles in the denser arrangement, and vice versa. In our simulations, we often set the particle arrangement based on a reasonable prescribed threshold of angular momentum loss. However, to completely eliminate such loss, one can employ a projection algorithm as suggested in~\cite{Lee2016}. We note that we do not pursue this idea further in this paper.} 

\begin{figure}[]
	\centering
	\includegraphics[width=0.7\textwidth]{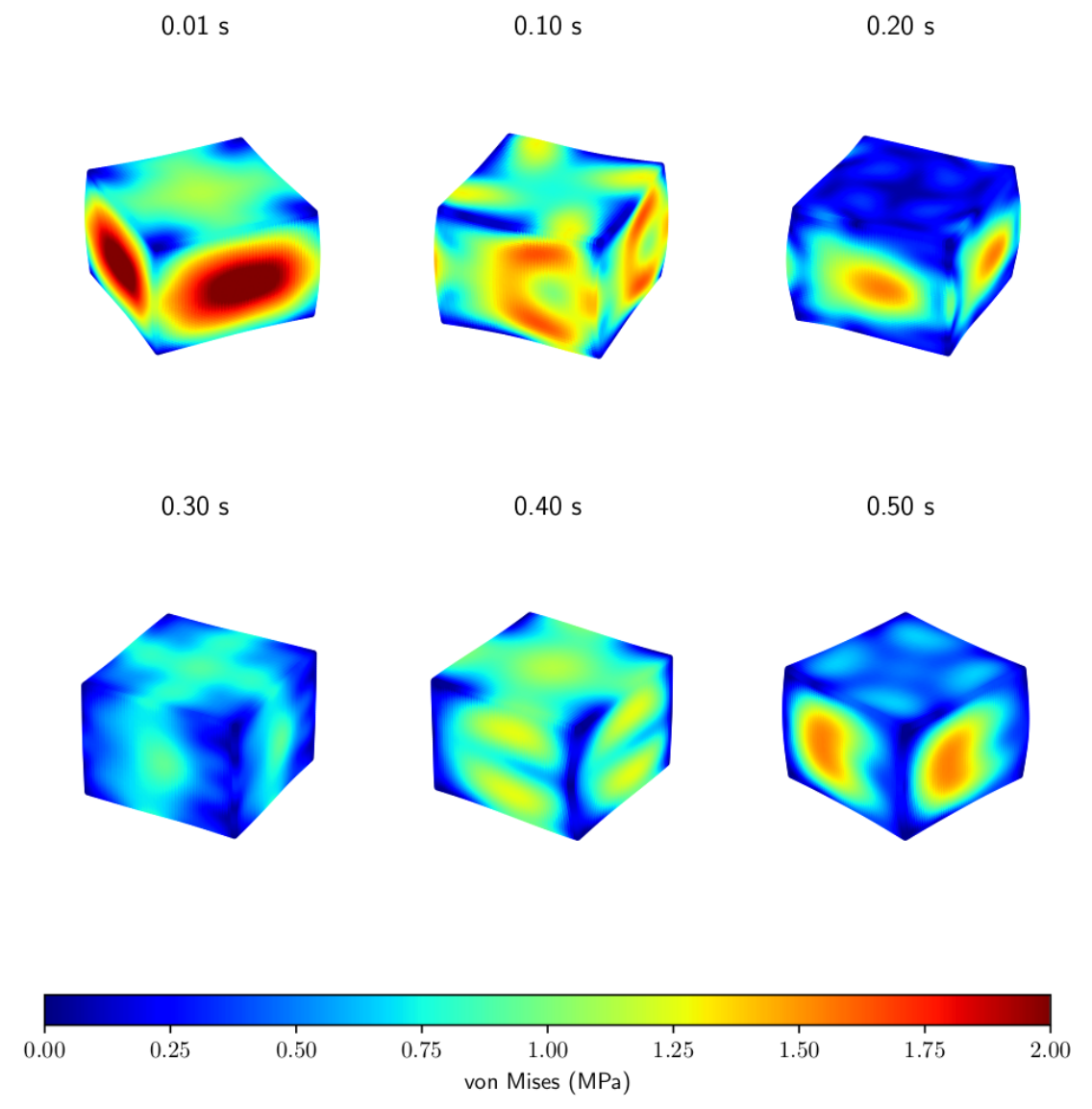}
	\caption{Stress contours of the spinning cube at various time frames}
	\label{fig:spinning_cube_contours}
\end{figure}

\begin{figure}[]
	\centering
	\includegraphics[width=1\textwidth]{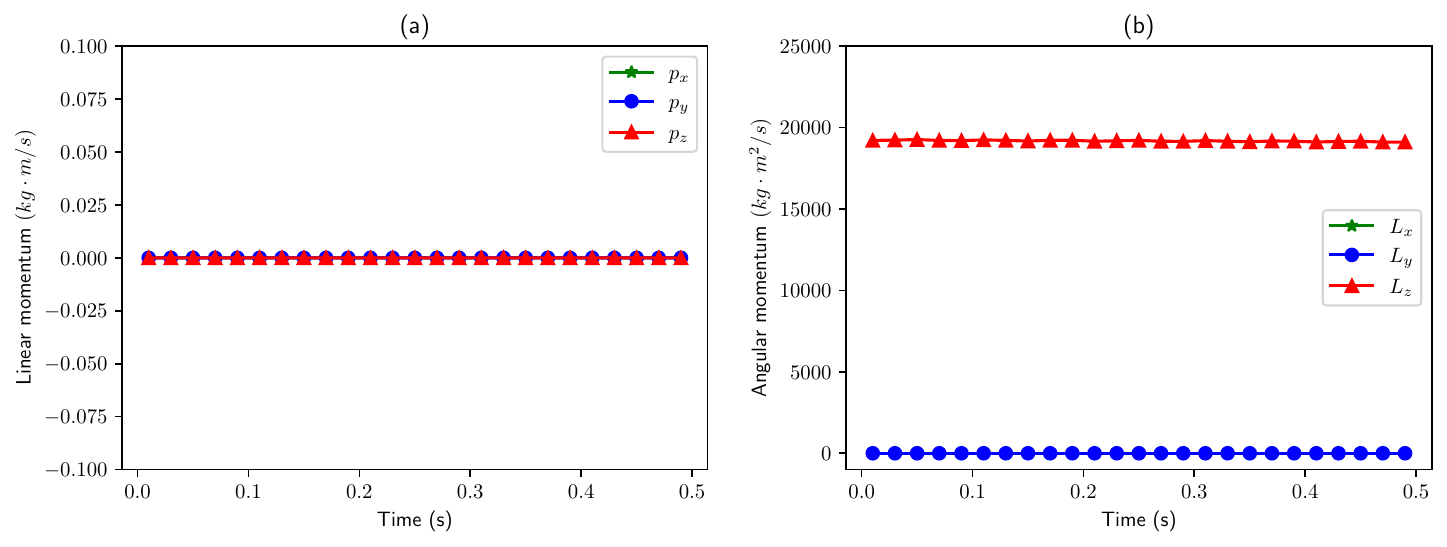}
	\caption{Linear and angular momentum of the spinning cube problem: (a) linear momentum, (b) angular momentum}
	\label{fig:spinning_cube_momentum}
\end{figure}

We extend the example to a three-dimensional rotating cube analogous to the previously discussed spinning plate. The cube has a width of 1 m and the material properties include a modulus of elasticity $E=17~\text{MPa}$, Poisson's ratio $\nu=0.3$, and density $\rho_o=1100~\text{kg/m}^3$, following the example shown in~\cite{Lee2016}. The cube is discretised using $50 \times 50 \times 50$ particles. We consider the same initial angular velocity as employed  in the 2D spinning plate problem. The von Mises stress distribution is presented in Figure~\ref{fig:spinning_cube_contours} from $t=0~\text{s}$ to $t=2~\text{s}$. {It is evident that the proposed method leads to the stress contour without spurious non-physical oscillation in 3D space, and the total linear and angular momentum of the system are shown in Figure~\ref{fig:spinning_cube_momentum}.} 

\subsection{Energy preservation analysis using finite deformation of a column}
\label{sec:column-3d}
\begin{figure}[ht!]
	\centering
	\includegraphics[width=0.4\textwidth]{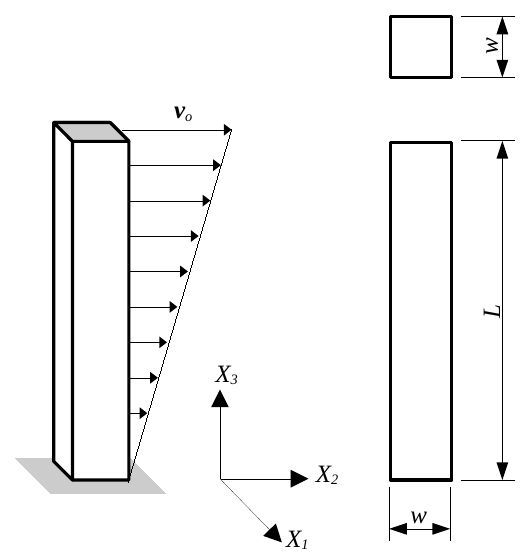}
	\caption{Configuration and dimension of the column for energy preservation analysis}
	\label{fig:cb_specs}
\end{figure}

In this example we assess the preservation of energy using a square column subjected to an initial velocity, as previously demonstrated in~\cite{Lee2016}. The dimensions of the column are $1 \times 1 \times 6~\text{m}$ and it is fixed at the bottom end, as depicted in Figure~\ref{fig:cb_specs}. 
The smoothing kernel used in this example is a cubic b-spline with smoothing length to particle spacing ratio $\beta = 0.9$. The column is made from a rubber material with modulus of elasticity $E=17~\text{MPa}$, {Poisson's ratio $\nu=0.495$}, and density $\rho_o=1100~\text{kg/m}^3$. An initial velocity $\vec{v}_o = (0, 5 X_3 /3, 0)^\trans ~\text{m/s}$ is applied to the beam at time $t=0~\text{s}$.

\begin{figure}[]
	\centering
	\includegraphics[width=1\textwidth]{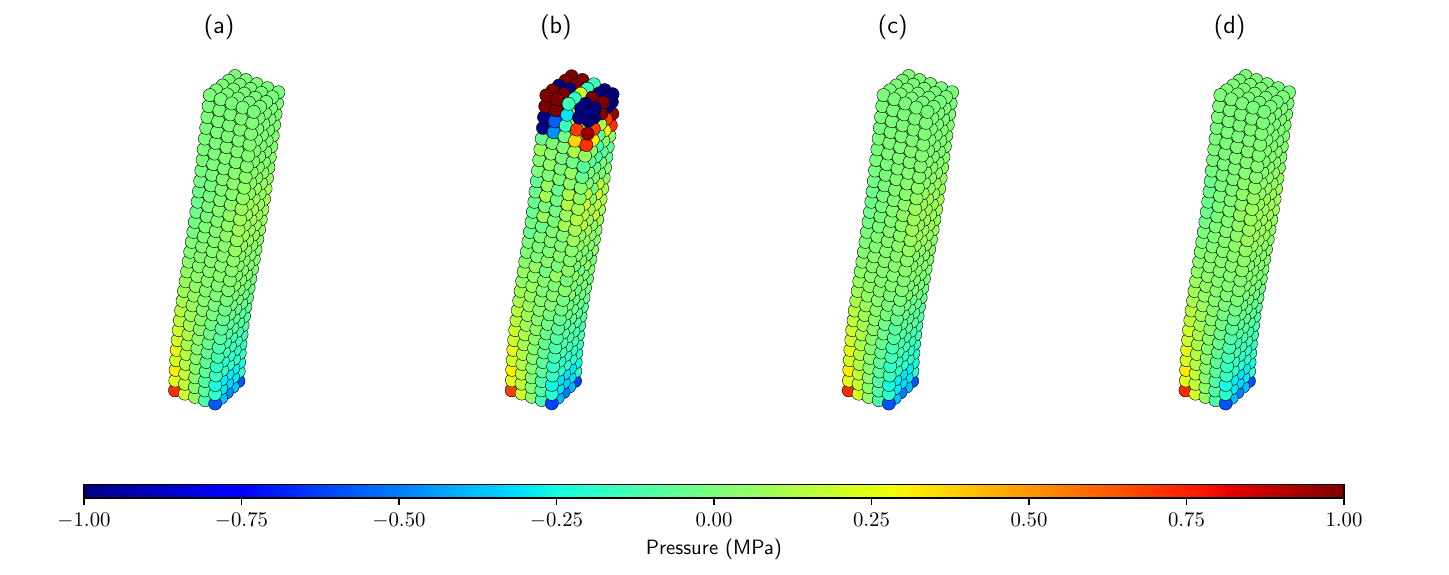}
	\caption{Deformation and pressure distribution at time $t$ = 75.05 ms using $5 \times 5 \times 30$ particle configuration, $\eta^{(2)}=0$ and $\eta^{(4)}=\frac{1}{64}$: (a) JST SPH, $C_{p,ij}$, $\alpha_{CFL}$ = 0.3, (b) JST SPH, $C_{p,ij}$, $\alpha_{CFL}$ = 0.5, (c) Bond-based SPH, $C_p^{Lin}$, $\alpha_{CFL}$ = 0.9, (d) Bond-based SPH, $C_{p,ij}$, $\alpha_{CFL}$ = 0.9}
	\label{fig:column_jst_fij}
\end{figure}

{We first compare the performance of our proposed bond-based SPH with JST-SPH approach described in~\cite{Lee2016} with respect to $\alpha_{CFL}$ and pressure wave speed $C_p$, as shown in Figure~\ref{fig:column_jst_fij}. In this example, bond-based SPH converges even with high $\alpha_{CFL}$=0.9, as seen in Figures~\ref{fig:column_jst_fij}(c) and (d) without encountering instabilities. On the other hand, JST-SPH experiences instabilities and shows spurious mode at $\alpha_{CFL}$=0.5 as illustrated in Figure~\ref{fig:column_jst_fij}(b). The spurious mode can be avoided by using lower $\alpha_{CFL}$=0.3 for JST-SPH, as shown in Figure~\ref{fig:column_jst_fij}(a). Furthermore, in this example, the bond-based SPH converges with $\alpha_{CFL}$=0.9 using $C_p$ with the linear term ($C_p^{Lin}$) alone (Figure~\ref{fig:column_jst_fij}(c)), as well as with both linear and nonlinear terms $C_{p,ij}$ (Figure~\ref{fig:column_jst_fij}(d)). In contrast, JST-SPH requires both linear and nonlinear terms to be active to produce results shown in Figures~\ref{fig:column_jst_fij}(a) and (b). The time increments in these simulations are presented in Table~\ref{table:cp_dt}. Here we use $C_p^{Lin}$ to compute the time increment for bond-based SPH since the differences in time increment between $C_p^{Lin}$ and $C_{p,ij}$ are not significant. Moreover, calculating the stretch at each particle to compute $C_{p,ij}$ also increases computational effort.
\begin{table}[ht!]
	\centering
	\begin{tabular}{|c|c|c|}
		\hline
		$\alpha$ & $C_p$ & $\Delta t$ ($\mu$s) \\
		\hline
		0.3 & $C_p^{Lin}$ & 73.75  \\
		0.3 & $C_{p,ij}$  & 83.04  \\
		0.9 & $C_p^{Lin}$ & 221.3  \\
		0.9 & $C_{p,ij}$  & 249.1  \\
		\hline
	\end{tabular}
	\caption{Magnitude of time increment due to variation of $C_p$ and $\alpha_{CFL}$ on $5 \times 5 \times 30$ particle configuration}
	\label{table:cp_dt}
\end{table}
}

\begin{figure}[]
	\centering
	\includegraphics[width=1\textwidth]{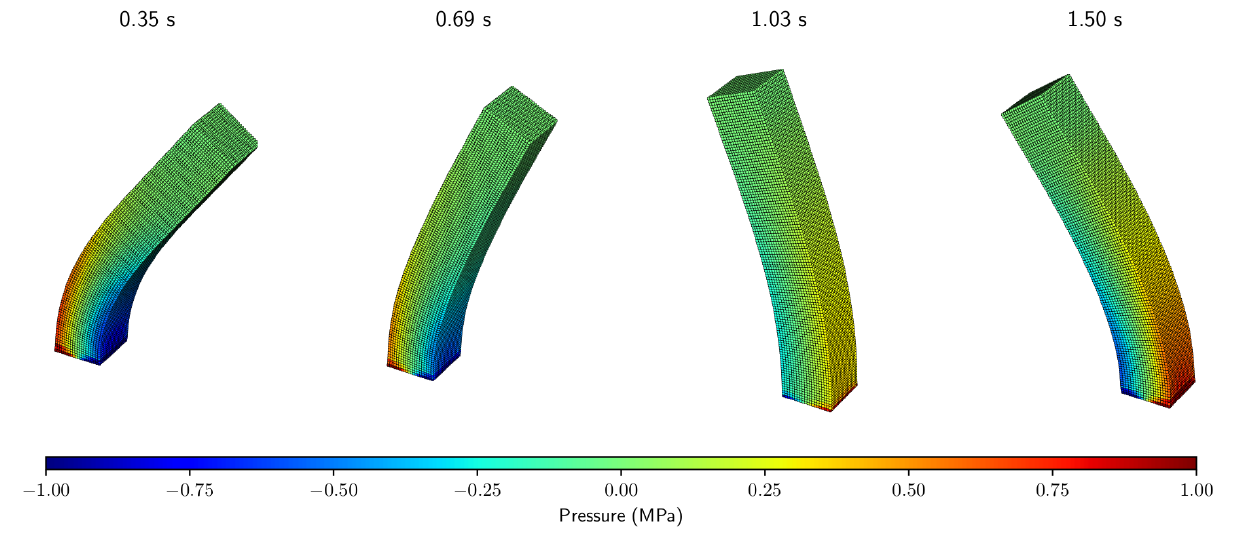}
	\caption{Deformation and pressure distribution of the column at various time frames using $20 \times 20 \times 120$ particle configuration, $\eta^{(2)}=0$ and $\eta^{(4)}=\frac{1}{64}$}
	\label{fig:bending_column}
\end{figure}

\begin{figure}[]
	\centering
	\includegraphics[width=1\textwidth]{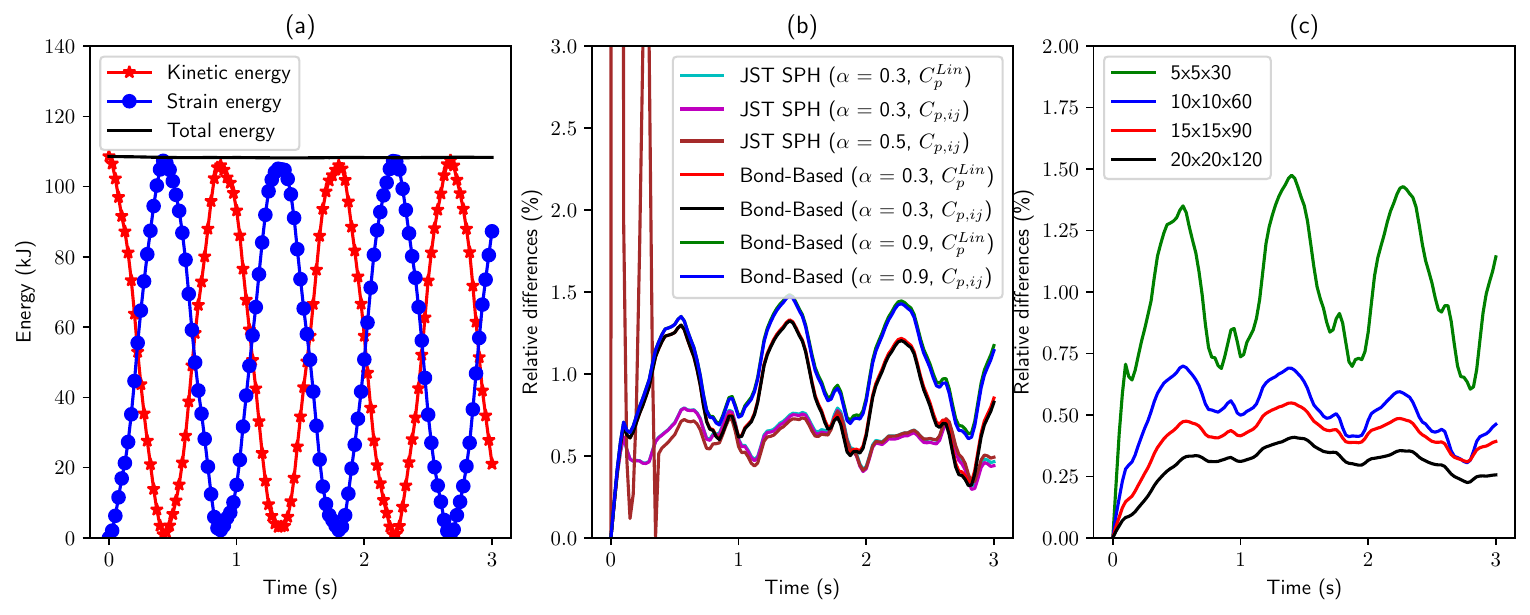}
	\caption{Energy analysis of the column problem: (a) energy history using Bond-based SPH with $20 \times 20 \times 120$ particles, (b) relative differences of total energy obtained from JST SPH and Bond-based SPH with $5 \times 5 \times 30$ particle configuration, (c) effect of particle refinement on energy preservation using Bond-based SPH}
	\label{fig:bending_energy}
\end{figure}

The total energy of the system can be decomposed into kinetic energy $E_k$ and strain energy $W_s$ where each can be obtained by summing up the individual particle contributions. The particle's kinetic energy is obtained by considering its mass and velocity. The particle's strain energy is determined based on the hyperelastic Neo-Hookean model, i.e., 
\begin{subequations}
	\begin{align}
		E_k &= \frac{1}{2} \sum_{i=1}^N  m_i \left(\vec{v}_i \cdot \vec{v}_i \right) \, \\
		W_s &= {\sum_{i=1}^N V_{o,i} \left( \frac{\mu}{2} \left( tr \, \widehat{\vec{C}}_i - 3 \right) + \frac{\kappa}{2} \left(J_i-1 \right)^2 \right); \hspace{1cm} \widehat{\vec{C}}_i = det\left( \vec{C}_i \right)^{-1/3} \vec{C}_i} \\
		E_t &= E_k + W_s \, ,
	\end{align}
	\label{eq:energy_function}
\end{subequations}
where $\mu$ stands for the shear modulus, $\kappa$ is the bulk modulus, $J$ is the Jacobian of the deformation gradient, and $\vec{C}$ is the right Cauchy-Green deformation tensor. 
The deformed shape and the {pressure} distribution of the column with $20 \times 20 \times 120$ particle configuration at various time frames are displayed in Figure~\ref{fig:bending_column}. It is evident that the contours are free from spurious oscillations and qualitatively show agreement with the reference~\cite{Lee2016}. {The changes in kinetic, strain, and total energy of the system are depicted in Figure ~\ref{fig:bending_energy}(a). At certain time instances, the peak kinetic energy does not reach the total energy value, indicating that some energy is stored in the form of strain energy. Similarly, there are occasions when the peak strain energy does not reach the total energy, because parts of the beam still remain in motion. This energy exchange often influenced by factors such as initial perturbation and dynamics characteristics of the column.}
	
{Figure ~\ref{fig:bending_energy}(b) illustrates the comparison of energy preservation properties between JST-SPH and bond-based SPH using a $5 \times 5 \times 30$ particle configuration. The relative difference of total energy at each time frame ($E_t^n$) compared to the initial total energy ($E_t^0$) is calculated as follows
\begin{equation}
	\text{Relative differences (\%)} = \frac{E_t^n - E_t^0}{E_t^0} \times 100 \, .
\end{equation}
It is evident from Figure ~\ref{fig:bending_energy}(b) that the energy loss is approximately~$1~\%$ for most cases, except for JST-SPH using $\alpha_{CFL}=0.5$ with $C_{p,ij}$ active. This total energy loss is attributed to the numerical dissipation induced by JST stabilisation terms. Furthermore, we study the effect of particle refinement on energy loss, as shown in Figure~\ref{fig:bending_energy}(c). It is evident that denser particle distributions result in lower total energy loss.}



\subsection{Stability analysis using pulling test of column}
\label{sec:pulling}
\begin{figure}[]
	\centering
	\includegraphics[width=0.4\textwidth]{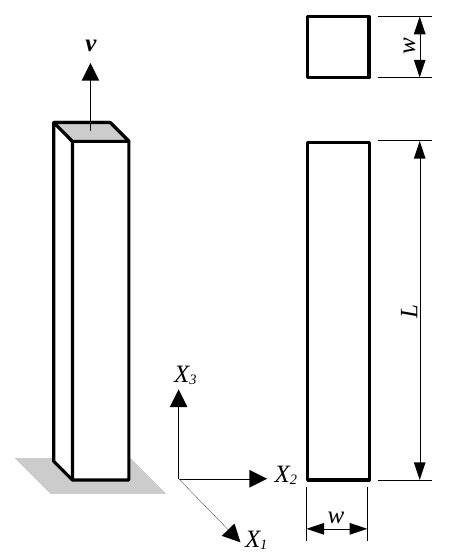}
	\caption{Configuration and dimension of the column for the pulling test}
	\label{fig:column_specs}
\end{figure}

In this example we assess the numerical stability of the proposed TLSPH method using a three-dimensional column subjected to extreme elongation. The square column has a width of $w = 1~\text{m}$ and has initial length of $L= 6~\text{m}$. The physical fields are discretised using $12 \times 12 \times 72$ uniformly distributed particles. Here we use a cubic b-spline kernel with first order correction and a width of $h = 0.5$.  The column is made from rubber with elasticity modulus $E=17~\text{MPa}$, and  $\rho=1100~\text{kg/m}^3$. The bottom end of the column is fixed and the free top end is subjected to a constant velocity $\vec{v} = (0,0,10)~\text{m/s}$, as shown in Figure~\ref{fig:column_specs}. 

{In this example, we vary Poisson's ratio to analyse the stability enhancement provided by bond-based SPH in handling extreme elongation due to tension loading. We consider two values of Poisson's ratio, $\nu= (0.45, \, 0.495)$, along with JST parameters, $\eta^{(2)}=0.0$ and $\eta^{(4)}=\frac{1}{64}$. The maximum stretch in $X_3$ direction, $\lambda_3$, is calculated as
\begin{equation}
	\lambda_3 = \frac{L' - L}{L} \times 100 \, ,
\end{equation}
where $L'$ is the final length of the column before any instability occurs. In our simulation this is indicated when the calculation of displacement and/or stress results in "NaN".
}

\begin{figure}[]
	\centering
	\includegraphics[width=0.5\textwidth]{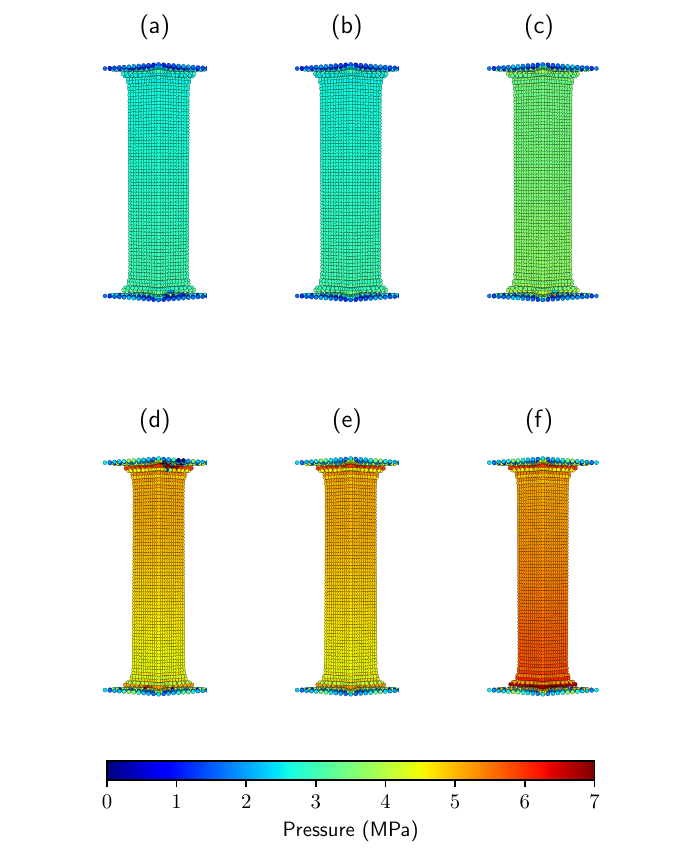}
	\caption{Deformation and pressure distribution of a column which pulled by a constant speed: (a) JST SPH, $\nu$=0.45, $\lambda_3$=421.90\%, (b) Bond-based SPH, $\nu$=0.45, $\lambda_3$=421.90\%, (c) Bond-based SPH, $\nu$=0.45, $\lambda_3$=526.80\%, (d) JST SPH, $\nu$=0.495, $\lambda_3$=427.58\%, (e) Bond-based SPH, $\nu$=0.495, $\lambda_3$=427.58\%, and (f) Bond-based SPH, $\nu$=0.495, $\lambda_3$=464.50\%}
	\label{fig:pulling_test}
\end{figure}

{We compare the maximum stretch of the column at which instabilities begin using both JST-SPH~\cite{Lee2016} and our proposed bond-based approach. For $\nu$=0.45, the particle deformations and the respective pressure distribution are illustrated in Figure~\ref{fig:pulling_test}~(a-c). Figure~\ref{fig:pulling_test}(a) shows the result of JST-SPH, where instability occurs when $\lambda_3=421.90 \%$. In contrast, simulations using SPH with the improved bond-based deformation gradient remain stable even at $\lambda_3=421.90 \%$, as shown in Figure~\ref{fig:pulling_test}(b), and can continue until $\lambda_3=526.80 \%$ before encountering instability, see Figure~\ref{fig:pulling_test}(c). This indicates that the proposed method achieves $24.86\%$ higher tensile elongation than JST-SPH. The simulation results with $\nu$=0.495 are presented in Figure~\ref{fig:pulling_test}~(d-f). JST-SPH exhibits instabilities when reaching $\lambda_3=427.58 \%$ stretch, while the proposed bond-based SPH remains stable until $\lambda_3=464.50 \%$, i.e., the bond-based SPH achieves $8.64 \%$ higher elongation than JST-SPH. We note that the simulations with higher Poisson's ratio ($\nu$=0.495) result in much higher stress compared to those using $\nu$=0.45 even when the elongations are comparable. This is because a higher Poisson's ratio leads to greater area reduction under the same load and initial configuration.}

\subsection{Stability analysis using twisting column problem}
\label{sec:twisting}
\begin{figure}[]
	\centering
	\includegraphics[width=0.4\textwidth]{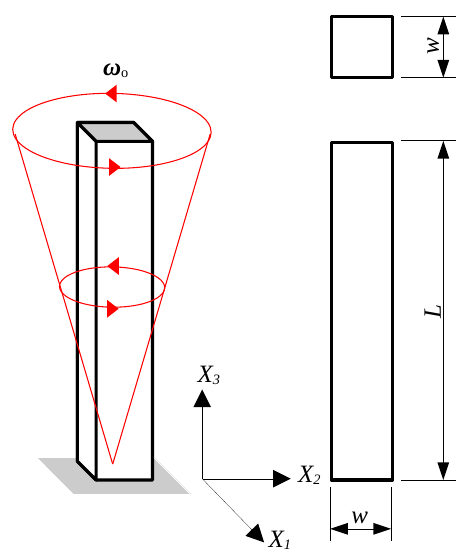}
	\caption{Configuration and dimension of the twisting column problem}
	\label{fig:twisting_specs}
\end{figure}

\begin{figure}[]
	\centering
	\includegraphics[width=1\textwidth]{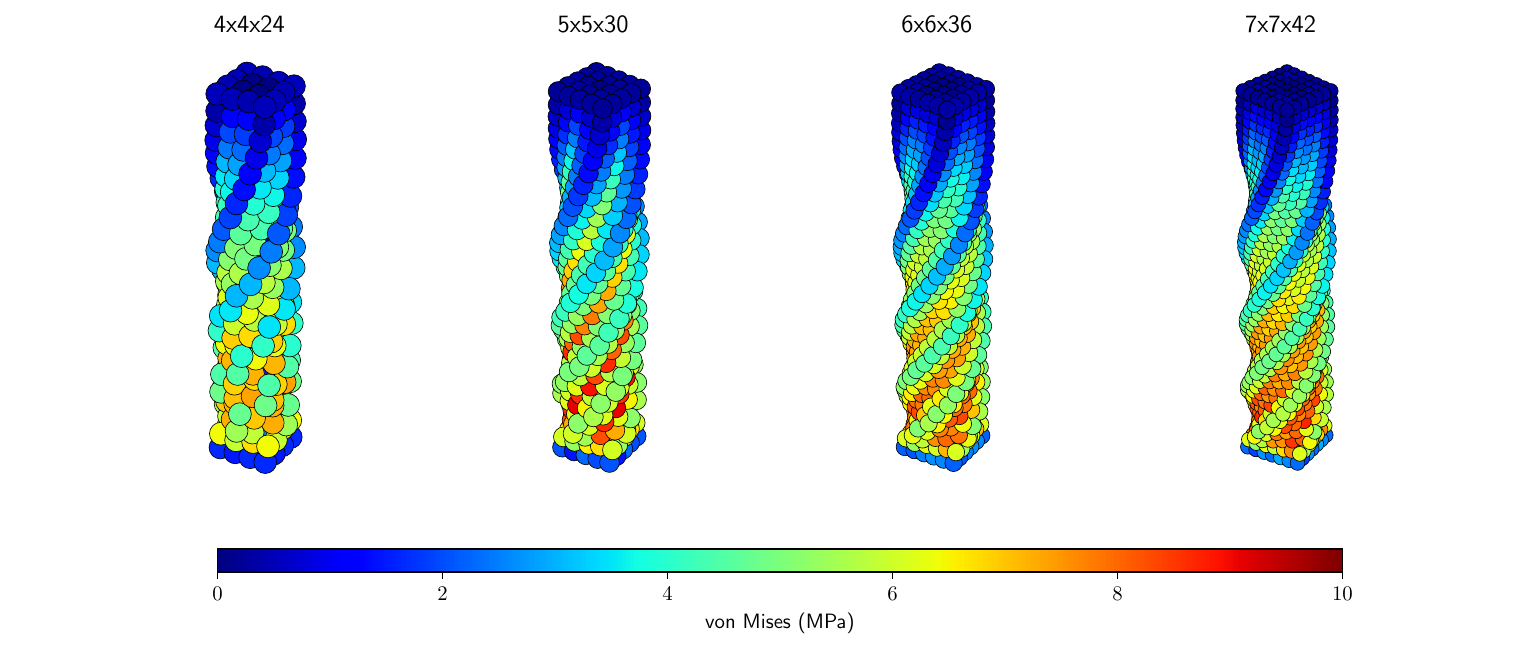}
	\caption{Analysis of particle refinement on the deformation of twisting column with $\omega_3 = 105$ $rad/s$, $E=17$ $MPa$, $\nu=0.45$, $\rho=1100$ $kg/m^3$, $\alpha_{CFL}=0.9$, $\eta^{(2)}=0$ and $\eta^{(4)}=0.125$. The results were recorded at $t=100$ ms}
	\label{fig:twisting_105_045}
\end{figure}

Our final example involves the twisting of a three-dimensional column with extremely large deformation. The square column is $1 \times 1 \times 6~\text{m}$ and is fixed at the bottom end, as depicted in Figure~\ref{fig:twisting_specs}. An initial angular velocity $\vec{\omega}_o=(0,0,\omega_3 \sin{\left( \pi X_3 / 12\right)})^\trans$ is applied to the column at time $t=0~\text{s}$. 
The deformation of the column and the von Mises stress for $\omega_3 = 105~\text{rad/s}$ are presented in Figure~\ref{fig:twisting_105_045} for various particle discretisations at time $t=100~\text{ms}$. It is evident that our proposed method exhibits the expected deformation patterns for all discretisations, even with a sparse $4 \times 4 \times 24$ particles. This is already an improvement with respect to the original TLSPH with JST stabilisation~\cite{Lee2016}, which requires a dense particle distribution to obtain a stable result with accurate deformation pattern. Our presented results are comparable with the upwind-SPH shown in~\cite{Lee2019}.

\begin{figure}[]
	\centering
	\includegraphics[width=1\textwidth]{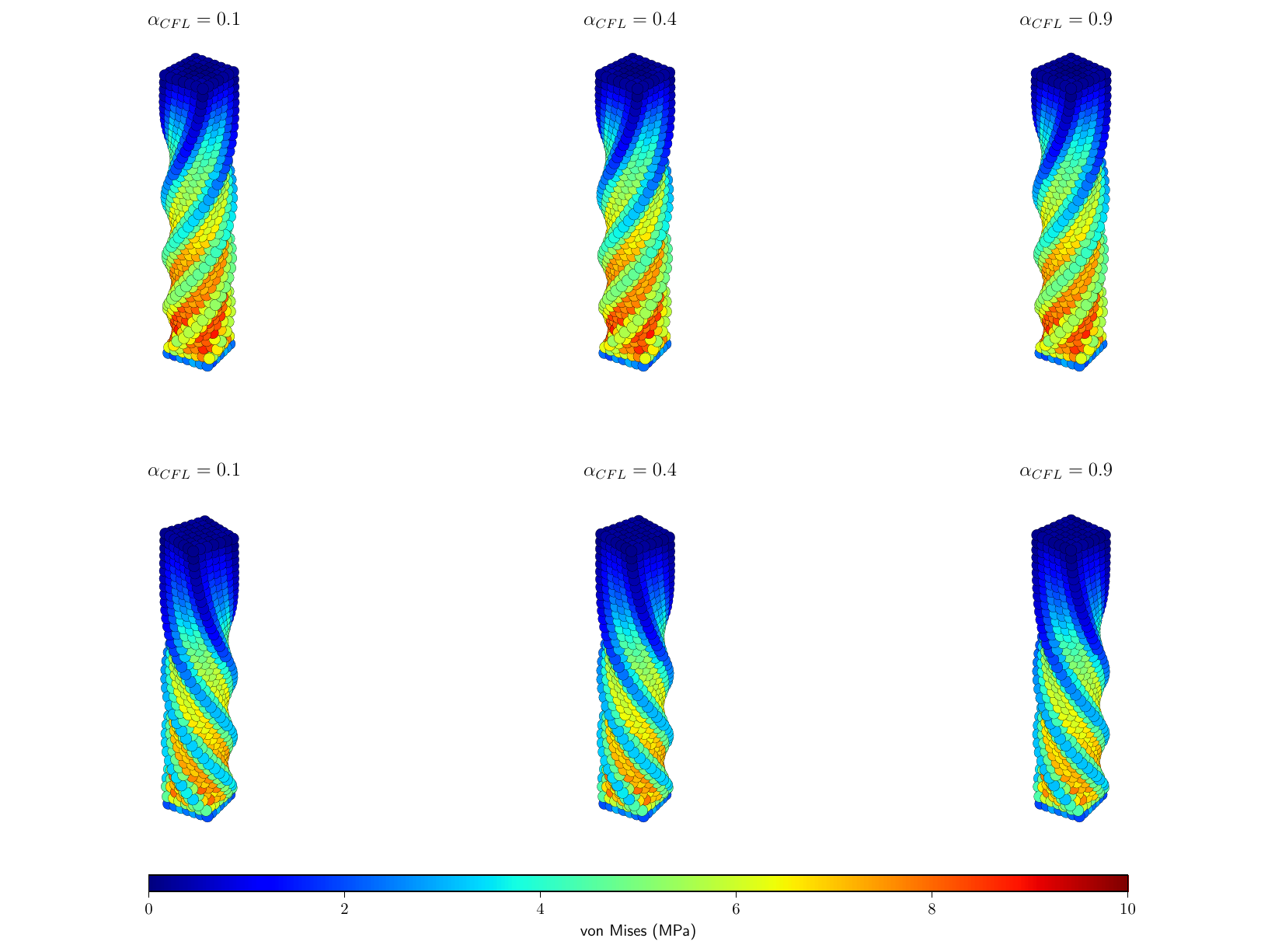}
	\caption{Analysis of time increment on the deformation of twisting column with particle configuration $7 \times 7 \times 42$, $\omega_3 = 105$ $rad/s$, $E=17$ $MPa$, $\nu=0.45$, $\rho=1100$ $kg/m^3$, $\eta^{(2)}=0$ and $\eta^{(4)}=0.125$. The upper row of the figure was recorded at $t=100$ ms, and the bottom row at $t=250$ ms}
	\label{fig:twisting_105_045_alpha}
\end{figure}

Next we investigate the effect of time increment on the numerical stability of the proposed method through the variation of the CFL constant $\alpha_{CFL} \in \{0.1,0.4,0.9\}$. The results are presented in Figure~\ref{fig:twisting_105_045_alpha} for $\omega_3 = 105~\text{m/s}$ recorded at time $t=100~\text{ms}$ (upper row) and $t=250~\text{ms}$ (bottom row). The results agree well with those presented in~\cite{Lee2016}. The use of the three-stage time integrator allows for the use of higher~$\alpha_{CFL}$, leading to a lower computational cost. Additionally, the numerical stability and accuracy are well-maintained, as evidenced by the consistent deformed shape and stress distribution.

\begin{figure}[ht!]
	\centering
	\includegraphics[width=1.0\textwidth]{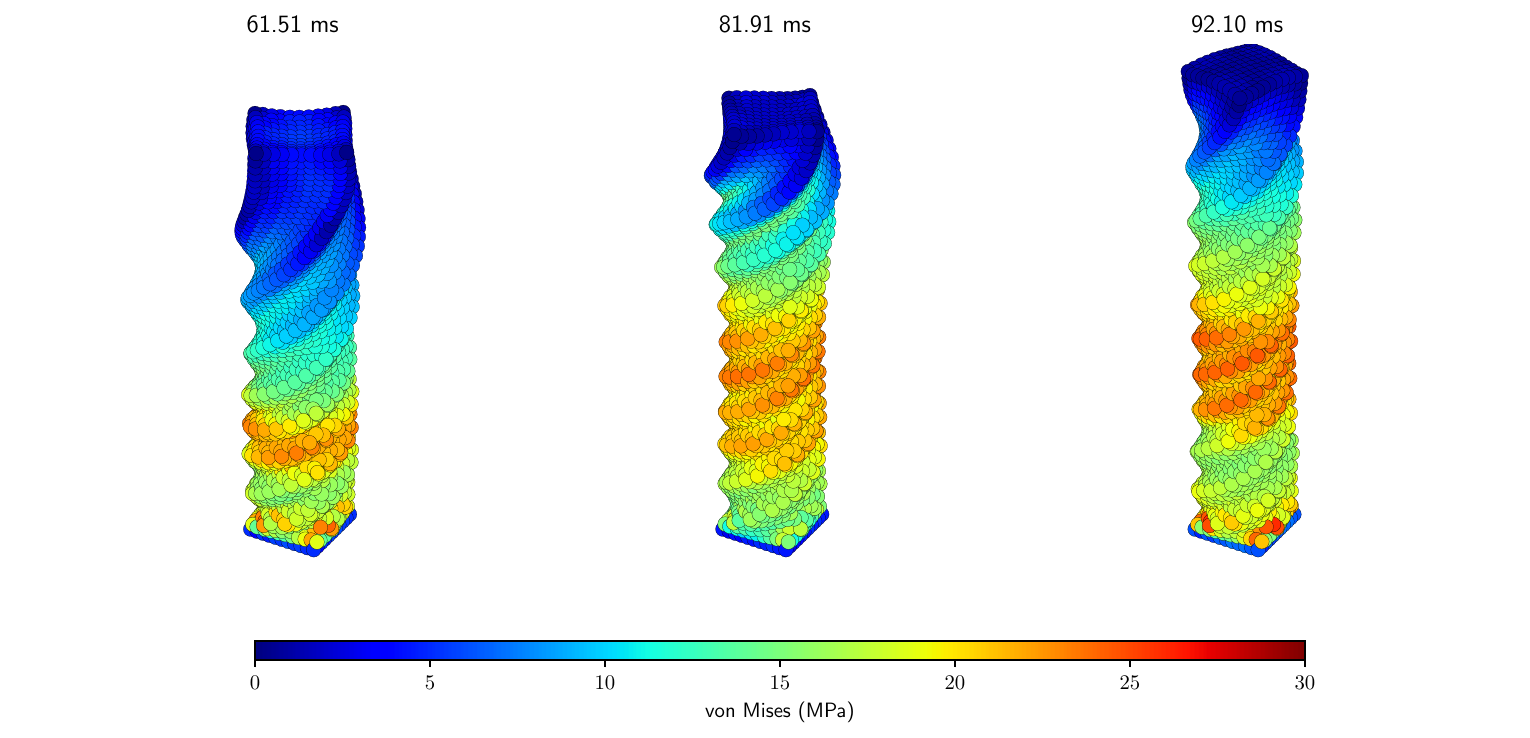}
	\caption{Deformation of twisting column with particle configuration $11 \times 11 \times 66$, $\omega_3 = 300$ $rad/s$, $E=17$ $MPa$, $\nu=0.49$, $\rho=1100$ $kg/m^3$, $\alpha_{CFL}=0.5$, $\eta^{(2)}=0$ and $\eta^{(4)}=0.125$.}
	\label{fig:twisting_300}
\end{figure}
\begin{figure}[ht!]
	\centering
	\includegraphics[width=1\textwidth]{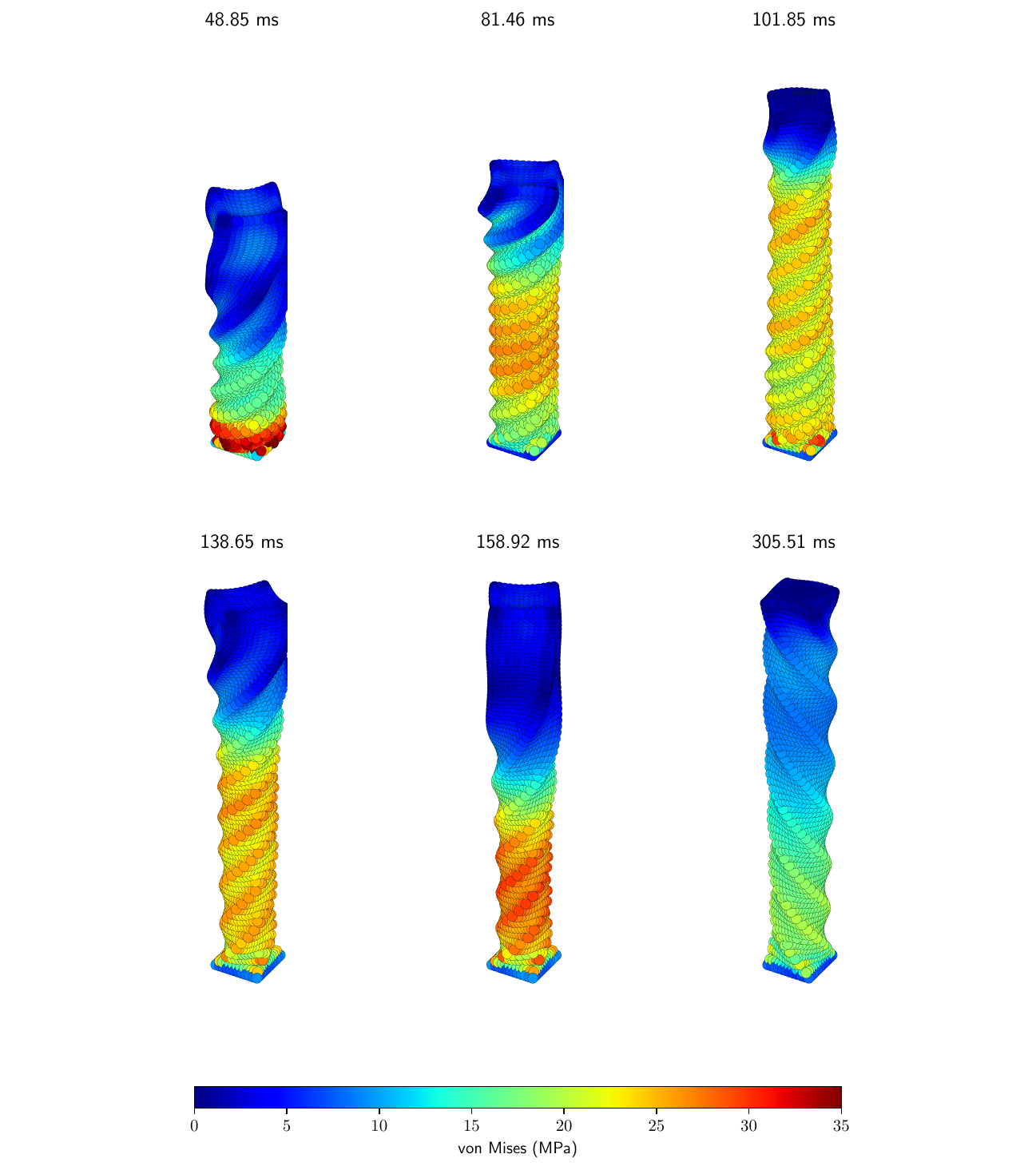}
	\caption{Deformation of twisting column with particle configuration $12 \times 12 \times 72$, $\omega_3 = 400$ $rad/s$, $E=17$ $MPa$, $\nu=0.49$, $\rho=1100$ $kg/m^3$, $\alpha_{CFL}=0.5$, $\eta^{(2)}=0$ and $\eta^{(4)}=0.6$.}
	\label{fig:twisting_400}
\end{figure}

Lastly we investigate the effect of the initial angular velocity applied to the column, denoted as $\omega_3$. Figure~\ref{fig:twisting_300} illustrates the deformed shape with $\omega_3=300~\text{rad/s}$. As reported in~\cite{Wu2023}, solving this problem with the standard SPH leads to numerical instability and particle detachment. In contrast, the proposed formulation exhibits stable results until $t = 92.1~\text{ms}$. To create an even more interesting and challenging test, we further increase the angular velocity to $\omega_3=400~\text{rad/s}$. The results are shown in Figure~\ref{fig:twisting_400}. It is evident that the simulation remains stable even for extremely large deformation. 

%% file: conclusions.tex
\section{Conclusions}
\label{sec:conclusions}

We presented a total Lagrangian Smoothed Particle Hydrodynamics (TLSPH) method designed for large-strain solid dynamics problems, showcasing optimal convergence and improved stability. Specifically, our approach involves a bond-based deformation gradient computed at particle pairs smoothed by a compactly supported Lagrangian kernel. Additionally, we consider the bond of a particle with itself to maintain the method's consistency established by the kernel correction. The concept was initially illustrated in one dimension and its excellent properties were demonstrated through the gradient estimation of polynomial and oscillatory deformation field. Intergation into the TLSPH framework involved incorporating it as a correction to the standard, or averaged, deformation gradient. Our proposed framework retains optimal convergence rates in $L_2$-norm and $H_1$-seminorm of displacement, stemming from the first order kernel correction, which ensures linear polynomial reproducibility. We have also shown that the proposed method preserves linear and angular momentum, as well as the total energy throughout simulations. Furthermore, our deformation gradient formulation exhibits increased stability in cases involving large strain, as illustrated in our pulling and twisting of column examples. In the pulling case, a strain up to $25\%$ larger can be achieved without spurious oscillations due to zero-energy modes, compared to TLSPH using the standard deformation gradient formulation. Notably, the proposed method requires the computation of the corrective bond-based gradient for each pair within the neighbourhood of a particle. However, due to its total Lagrangian construction, this computation can be efficiently performed without reidentifying the neighbours.

There are several promising future developments of the proposed method. The first concerns integrating the bond-based deformation gradient in the updated-Lagrangian SPH framework. This will enable computation of complex problems involving topological changes, such as crack propagation and free surface flow. {Another important aspect is the investigation of methodologies to ensure conservation of energy and angular momentum. As demonstrated in the examples, a small percentage of losses are observed due to dissipation and non co-linearity between pairwise force and the relative current position vector.} Another promising direction follows the findings in our one-dimensional analysis, which indicates that the proposed improved bond-based deformation gradient is more accurate for higher order polynomial or highly oscillatory fields. This has potential applications in unsteady analyses of both solid and fluid dynamics.

%% file: appendix.tex
\section{Gradient of the first order corrected kernel}
\label{sec:appendix}
In this appendix we present the derivation to obtain the gradient of the first order corrected kernel $\nabla \widetilde{W}^1_i (\vec x_j)$ as previously described in Section~\ref{sec:corrected-sph}. First we define a compact notation $\vec{x}_{ij} = \vec{x}_i - \vec{x}_j$ for the vector between two coordinates $\vec x_i$ and $\vec x_j$. We then rewrite equation~\ref{eq:cw} as
\begin{equation}
	\label{eq:cw1}
	\widetilde{W}_{i}^1 (\vec x_j) = \alpha_i \left( 1 +  \vec{\beta}_i \cdot \vec x_{ij} \right) W_i(\vec x_j) \, .
\end{equation}
The coefficients $\alpha_i$ and $\beta_i$ in~\eqref{eq:cw1} are described as
\begin{subequations}
	\begin{align}
		\alpha_i &= \left( \varphi_i - \left( (\vec{\Phi}_i)^{-1} \vec{\phi}_i \right) \cdot \vec{\phi}_i \right)^{-1} \, , \\
		\vec{\beta}_i &= - \left( \vec{\Phi}_i \right)^{-1} \vec{\phi}_i \, .
	\end{align}
\end{subequations}
where $\varphi_i$, $\vec{\phi}_i$, and $\vec{\Phi}_i$ are the geometric moments of the kernel centred at $\vec x_i$, which are given by
\begin{subequations}
	\begin{align}
		\varphi_i &= \sum_{j \in \set N(\vec x_i)} V_j W_i (\vec x_j) \, \\
		\vec{\phi}_i &= \sum_{j \in \set N(\vec x_i)} \vec{x}_{ij} V_j W_i (\vec x_j) \, \\
		\vec{\Phi}_i &= \sum_{j \in \set N(\vec x_i)} \vec{x}_{ij} \otimes \vec{x}_{ij} V_j W_i (\vec x_j) \, .
	\end{align}
\end{subequations}
The gradient of the first order corrected kernel $\nabla \widetilde{W}^1_i (\vec x_j)$ can be obtained by differentiating~\ref{eq:cw1}, which yields
\begin{equation}
	\label{eq:dcw1}
	\nabla \widetilde{W}_{i}^1 (\vec x_j) = 
	\alpha_i \left( 1 + \vec{\beta}_i \cdot \vec{x}_{ij} \right) \nabla {W}_{i} (\vec x_j) + 
	\nabla \alpha_i \left( 1 + \vec{\beta}_i \cdot \vec{x}_{ij} \right) {W}_{i} (\vec x_j) + 
	\alpha_i \left( (\nabla \vec{\beta}_i)^{\trans} \cdot \vec{x}_{ij} + \vec{\beta}_i \right) {W}_{i} (\vec x_j) \, .
\end{equation}
The derivative of $\alpha_i$ and $\vec{\beta}_i$ with respect to coordinate axis $k$ are given as
\begin{subequations}
	\begin{align}
	\label{eq:dalpha}
	{\nabla}_k \alpha_i &= - \alpha_i^2 \left( {\nabla}_k \varphi_i - \left((\vec{\Phi}_i)^{-1} \vec{\phi}_i \right) \cdot {\nabla}_k \vec{\phi}_i - \left((\vec{\Phi}_i)^{-1} {\nabla}_k \vec{\phi}_i \right) \cdot \vec{\phi}_i + \left((\vec{\Phi}_i)^{-1} {\nabla}_k \vec{\Phi}_i (\vec{\Phi}_i)^{-1} \vec{\phi}_i \right) \cdot \vec{\phi}_i \right) \, , \\
	\label{eq:dbeta}
	{\nabla}_k \vec{\beta}_i &= - (\vec{\Phi}_i)^{-1} {\nabla}_k \vec{\phi}_i + (\vec{\Phi}_i)^{-1} {\nabla}_k \vec{\Phi}_i (\vec{\Phi}_i)^{-1} \vec{\phi}_i \, .
	\end{align}
\end{subequations}
In the above, the derivate of the geometric moments are given as follows
\begin{subequations}
	\begin{align}
		\vec{\nabla} \varphi_i &= \sum_{j \in \set N(\vec x_i)} V_j \vec{\nabla} W_i (\vec x_j) \\
		\vec{\nabla} \vec{\phi}_i &= \sum_{j \in \set N(\vec x_i)} V_j \left( \vec{x}_{ij} \otimes \vec{\nabla} W_i (\vec x_j) + \vec{I} W_i (\vec x_j) \right) \\
		{\nabla}_k \vec{\Phi}_i &= \sum_{j \in \set N(\vec x_i)} V_j \left( \vec{x}_{ij} \otimes \vec{x}_{ij} {\nabla}_k W_i (\vec x_j) + \left(\vec{x}_{ij} \otimes \vec{\delta}^k + \vec{\delta}^k \otimes \vec{x}_{ij} \right) W_i (\vec x_j) \right) \, ,
	\end{align}
\end{subequations}
where $\vec{I}$ and $\vec{\delta}^k$ denote an identity matrix and a unit vector in the $k-$th direction, respectively. Note that although an explicit formulation can be obtained, the computation of the first derivative $\nabla \widetilde{W}^1_i (\vec x_j)$ is more involved. This includes obtaining the tensor ${\nabla}_k \vec{\Phi}_i $ and a small matrix inversion $(\vec{\Phi}_i)^{-1}$ for each position $\vec x_i$. 


\section{Derivation of the symmetrised conservation of linear momentum}
\label{sec:derivation_com}

{
Following the reference~\cite{frontiere2017}, we integrate the conservation of linear momentum over a control volume $V$ without external forces. Let us define a kernel interpolant $\psi_j = V_{o,j} \widetilde{W}_j^1(X_i)$ such that
\begin{equation}
    \label{eq:b1}
	\int_{V_o} \psi \frac{\D \vec{p}}{\D t} = \int_{V_o} \psi \, \text{div}(\vec{P}) \, .
\end{equation}
Here we employ a one-point quadrature approximation or the mid-point rule, that is $\int f \psi \approx V_{o,i} f_i$ for an arbitrary function $f$. Hence, the integral evaluation of~\eqref{eq:b1} over a small control volume at \textit{i} becomes
\begin{equation}
	\label{eq:integration-at-i}
	V_{o,i} \frac{\D \vec{p}_i}{\D t} =  \sum_{j \in \set N(\vec X_i)} \vec{P}_j \int_{V_o} \psi_i \nabla_o \psi_j \, .
\end{equation}
Integration in~\eqref{eq:integration-at-i} can be conducted by performing integration by parts, resulting in
\begin{equation}
	\label{eq:integration-bypart}
	V_{o,i} \frac{\D \vec{p}_i}{\D t} =  \sum_{j \in \set N(\vec X_i)} \vec{P}_j \left( \psi_i \int_{V} \nabla_o \psi_j - \int_{V_o} \psi_j \nabla_o \psi_i \right) \, .
\end{equation}
Then, divergence theorem is applied to ~\eqref{eq:integration-bypart}, which allows it to be rewritten as
\begin{equation}
	\label{eq:using-divergence}
	V_{o,i} \frac{\D \vec{p}_i}{\D t} =  \sum_{j \in \set N(\vec X_i)} \vec{P}_j \left( \oint_{\partial V_o} \psi_i \psi_j \vec{n} - \int_{V_o} \psi_j \nabla_o \psi_i \right) \, .
\end{equation}
Recall the kernel identity to modify the last term in ~\eqref{eq:using-divergence}, resulting in
\begin{equation}
	\label{eq:adding-Pi}
	\sum_{j \in \set N(\vec X_i)} \int_{V_o} \psi_i \nabla_o \psi_j = 0 \longrightarrow \sum_{j \in \set N(\vec X_i)} \vec{P}_i \int_{V_o} \psi_i \nabla_o \psi_j = 0 \, .
\end{equation}
Then, integration by parts and the divergence theorem are employed to derive the following relation
\begin{equation}
	\label{eq:adding-Pi-divergence}
	\sum_{j \in \set N(\vec X_i)} \vec{P}_i \int_{V_o} \psi_i \nabla_o \psi_j = 0 \longrightarrow \sum_{j \in \set N(\vec X_i)} \vec{P}_i \left( \oint_{\partial V_o} \psi_i \psi_j \vec{n} - \int_{V_o} \psi_j \nabla_o \psi_i \right) = 0 \, .
\end{equation}
Summation of Equation ~\eqref{eq:integration-at-i}, ~\eqref{eq:using-divergence}, ~\eqref{eq:adding-Pi}, and ~\eqref{eq:adding-Pi-divergence} leads to
\begin{equation}
	\label{eq:summationPiPj}
	\begin{aligned}
		2 V_{o,i} \frac{\D \vec{p}_i}{\D t} =  
		&\sum_{j \in \set N(\vec X_i)} \vec{P}_j \int_{V_o} \psi_i \nabla_o \psi_j + \sum_{j \in \set N(\vec X_i)} \vec{P}_j \left( \oint_{\partial V_o} \psi_i \psi_j \vec{n} - \int_{V_o} \psi_j \nabla_o \psi_i \right) + \\
		& \sum_{j \in \set N(\vec X_i)} \vec{P}_i \int_{V_o} \psi_i \nabla_o \psi_j + \sum_{j \in \set N(\vec X_i)} \vec{P}_i \left( \oint_{\partial V_o} \psi_i \psi_j \vec{n} - \int_{V_o} \psi_j \nabla_o \psi_i \right) \, .
	\end{aligned}
	\end{equation}
Equation ~\eqref{eq:summationPiPj} can be rearrange into
\begin{equation}
	\label{eq:symmetrisedPiPj_surface}
	2 V_{o,i} \frac{\D \vec{p}_i}{\D t} = \sum_{j \in \set N(\vec X_i)} \left( \vec{P}_i + \vec{P}_j \right) \left( \oint_{\partial V_o} \psi_i \psi_j \vec{n} + \int_{V} \psi_j \nabla_o \psi_i - \int_{V} \psi_i \nabla_o \psi_j \right) \, .
\end{equation}
Following the consistency condition, the surface integration on the right-hand side will naturally vanish.
\begin{equation}
	\label{eq:mid-point-rule-derivative}
	\int \psi_i \nabla_o \psi_j \approx V_{o,i} \nabla_o \psi_j \, .
\end{equation}
Recall that $\psi_j = V_j W_j(X_i)$ and one-point quadrature approximation in Equation ~\eqref{eq:mid-point-rule-derivative}, then, Equation ~\eqref{eq:symmetrisedPiPj} can be expressed as
\begin{equation}
	\label{eq:symmetrisedPiPj}
	2 V_{o,i} \frac{\D \vec{p}_i}{\D t} = \sum_{j \in \set N(\vec X_i)} \left( \vec{P}_i + \vec{P}_j \right) \left(V_{o,j} V_{o,i} \nabla \widetilde{W}_j^1(X_i) - V_{o,i} V_{o,j} \nabla \widetilde{W}_i^1(X_j) \right) \, .
\end{equation}
Finally, the conservation of linear momentum can be written as follow
\begin{equation}
	\label{eq:symmetrised_com}
	\frac{\D \vec{p}_i}{\D t} = \frac{1}{2} \sum_{j \in \set N(\vec X_i)} V_{o,j} \left( \vec{P}_i + \vec{P}_j \right) \left(\nabla_o \widetilde{W}_j^1(X_i) - \nabla_o \widetilde{W}_i^1(X_j) \right) \, .
\end{equation}
}